\documentstyle[11pt,oldlfont]{article}
\input{mac}
\input{serg.tex}

\begin{document}

\vskip .2in

\centerline{\bf The STRUCTURE of BRANCHING}
\centerline{\bf in ANOSOV FLOWS of $3$-MANIFOLDS}

\vskip .1in

\centerline{S\'ergio R. Fenley
\footnote{Research supported by NSF grants
DMS-9201744 and DMS-9306059.}}
\centerline{Mathematical Sciences Research Institute and}
\centerline{University of California, Berkeley}

\vskip .2in



\section{Introduction}

In this article we study the topological structure  of the lifts to
the universal of the stable and unstable foliations 
of $3$-dimensional Anosov flows.
In particular we consider the case when these foliations do not have
Hausdorff leaf space. 
We completely determine the structure of the set of 
non separated leaves from a given leaf in one of these foliations.
As a consequence of this suspensions are characterized, 
up to topological conjugacy,
as the only $3$-dimensional Anosov
flows without freely homotopic closed orbits.
Furthermore the structure of branching is
related to the topology of the manifold: 
if there are infinitely many leaves
not separated from each other, 
then there is an incompressible
torus transverse to the flow. 
Transitivity is not assumed for these results.
Finally,  if the manifold has negatively curved fundamental group 
we derive some important properties of the limit sets of leaves
in the universal cover.

\vskip .1in

This article deals with a powerful technique for analysing  
Anosov flows in dimension $3$, namely the study of 
the topological structure of the (weak) stable and
unstable foliations when lifted to the universal cover.
This technique was introduced in a remarkable paper of
Verjovsky \cite{Ve} in order to study codimension one
Anosov flows. 
If the lifted (say) stable foliation has
Hausdorff leaf space, then it is homeomorphic to the
set of real numbers and we say that the stable foliation in
the manifold is $\rrrr$-covered. When both foliations
are $\rrrr$-covered the flow is said to be 
$\rrrr$-covered. Two early uses of this technique were:
(1) Ghys \cite{Gh} proved that an Anosov flow in a Seifert
fibered space is $\rrrr$-covered. This was an essential
step in showing that the flow is, up to finite covers,
topologically conjugate to a geodesic flow in the unit
tangent bundle of a closed surface of negative curvature (briefly,
a geodesic flow). (2) If the fundamental group of the 
manifold is solvable then the $\rrrr$-covered property,
proved by Barbot \cite{Ba1,Ba2}, is again an essential step
in Plante's proof \cite{Pl2,Pl3} that the flow is
topologically conjugated to a suspension of an Anosov
diffeomorphism of the torus (a suspension). 
In fact this last result
holds for any codimension one Anosov flow.

More recently, a lot of information has been gained by
analysing not just the 
individual leaf spaces, but rather the joint
topological structure of the stable and unstable foliations.
Using this and Dehn
surgery on closed orbits of suspensions
or geodesic flows
\cite{Fr,Go}, 
a large family of examples
was constructed where every closed orbit of the flow
is freely homotopic to infinitely many other closed
orbits \cite{Fe3}. This never 
happens for suspensions or geodesic flows, and
was thought to be impossible for any Anosov flow.

In addition the topological study gives information about
metric properties
of flow lines: We say that a flow is
{\em quasigeodesic} if flow lines
are uniformly efficient (up to a bounded 
multiplicative distortion) in measuring distances
in relative homotopy classes. 
Suspensions and geodesic flows are always
quasigeodesic and there are
many quasigeodesic ``pseudo-Anosov" flows  in hyperbolic
$3$-manifolds \cite{Ca-Th,Mos}.
The Dehn surgery construction mentioned above produces
a large family of Anosov flows in hyperbolic
manifolds which are not quasigeodesic.

Barbot \cite{Ba3,Ba4} also used this topological theory
to study Anosov flows 
and proved the following remarkable result: Assume 
that there is a Seifert fibered piece of the torus decomposition of the  
manifold \cite{Jo,Ja-Sh} and suppose that the corresponding fiber is not
freely homotopic to a closed orbit of the flow. 
First isotopically adjust the boundary tori to be as transverse
to the flow as possible \cite{Ba3}. Then the 
flow in that piece is topologically conjugate to a (generalized) 
geodesic flow on the unit tangent bundle of a compact
surface with boundary.
If the manifold is a graph manifold 
and all fibers satisfy the condition above, then the flow is
up to topological conjugacy obtained by Dehn surgery
on finitely many closed orbits of a geodesic flow \cite{Ba5}.
Using this Barbot \cite{Ba4,Ba5}
has obtained the first known examples
of graph manifolds which are neither torus bundles
over the circle, nor Seifert fibered spaces
and which do not admit Anosov flows.

The results above are in great part due to a complete
characterization of the possible joint topological structures of
$\rrrr$-covered Anosov flows \cite{Ba2,Fe3}. On the other hand
very little is known about the non $\rrrr$-covered case, for the
simple reason that their structure is not understood
at all.
The purpose of this article is to start a systematic study of Anosov 
flows which are not $\rrrr$-covered, where we then say the lifted
foliations have branching.

It is easy to show that intransitivity implies that the 
flow is not $\rrrr$-covered
\cite{So,Ba1} and for many years there was a great effort in trying
to prove that these two properties are equivalent
\cite{Ve,Gh,Fe3,Ba2}. However in
a surprising development Bonatti-Langevin \cite{Bo-La} have
recently constructed a transitive, non $\rrrr$-covered Anosov 
flow in dimension $3$. Their example has an embedded torus
transverse to the flow.

This leads us to two basic and very important questions concerning 
branching: (1) when can branching occur and (2) what are the possible
structures of branching in Anosov flows of
$3$-manifolds. In this article we give a complete answer to the
second question. 
We then show that the structure of branching is strongly related to
dynamics of the flow, the topology of the manifold and the metric
behavior of the stable and unstable foliations.

Let then $\Phi$ be an Anosov flow in $M^3$ with two
dimensional stable and unstable foliations $\fs, \fu$.
Let $\fns, \fnu$ be the respective lifts to the
universal cover $\mi$.
Let $\hs$ and $\hu$ denote the leaf spaces of $\fns$ and
$\fnu$ respectively. If $\fs$ is not $\rrrr$-covered,
then $\hs$ is not Hausdorff. The {\em branching}
leaves of $\fns$ correspond to the non Hausdorff
points in $\hs$. Two leaves $F \not = F'$ of $\fns$
form a {\em branching pair} if the corresponding
points in $\hs$ are not separated from each other.
This is equivalent to saying that $F, F'$ do not
have disjoint saturated neighborhoods in $\mi$, where
a saturated neighborhood of $\fns$ is an open set
which is a union of leaves of $\fns$. Similarly for $\fnu$.

Since the universal cover is simply connected,
$\fns$ and $\fnu$ are always transversely orientable 
and an orientation is chosen.
Then there is a notion of branching in the positive
or negative directions. The first non trivial result
about the structure of branching is the following
\cite{Fe5}: Suppose the flow is transitive. If
there is branching in the positive direction
of (say) the stable foliation
then this foliation also has 
branching in the negative direction.
This is a result about the ``global" structure of branching.

We analyse the ``local" structure of branching.
For general foliations the branching of 
the lifted foliations can be very complicated \cite{Im}.
We will show here that branching in Anosov foliations is 
of a simple type which is very rigid.
For simplicity many theorems are stated for $\fns$ but they
work equally well for $\fnu$.

A leaf of $\fns$ or $\fnu$ is said to be
{\em periodic} if it is left invariant by  
a non trivial covering translation
of the universal cover.
Equivalently, its image in $M$ contains a closed
orbit of $\Phi$.
We first show that branching puts 
a restriction in the type of the leaf.

\vskip .1in
\noindent {\bf Theorem A} \
Let $\Phi$ be an Anosov flow in $M^3$. If $F$ 
is a branching leaf of $\fns$,
then $F$ is 
periodic. 
\vskip .1in

Theorem $A$ 
should be interpreted as a rigidity result in the sense
that periodic leaves are ``rigid" while non periodic leaves are non 
rigid. This is best seen in the manifold: if the stable
leaf (in the manifold) is periodic, then it contains
a closed orbit of $\Phi$ and every orbit in
the leaf is forward asymptotic
to this closed orbit. The nearby returns are in
the same local stable leaf. In case the leaf is
not periodic the forward orbits limit
in orbits in the manifold, 
but the nearby returns are always in distinct local stable leaves.
This means that when lifted to the universal
cover one can perturb slightly the local structure,
which will then produce a contradiction.

Our next goal is to understand the structure of the set ${\cal E}$ 
of non separated leaves from a given leaf $F$ of (say) $\fns$.
There is a natural order in ${\cal E}$ given by:
if $E, L \in {\cal E}$ then we say that $E < L$
if there are $G, H \in \fnu$ with $G \cap E \not =
\emptyset, H \cap L \not = \emptyset$ and $G$ is
in the back of $H$, see fig. \ref{bt}.
It is easy to see that this is a total order in ${\cal E}$.
Using this we can say that a branching leaf $D$ is
between $E$ and $L$ if $E < D < L$.


\blankfig{bt}{1.4}{The set of non separated leaves from $F \in \fns$.
$D$ is between $E$ and $L$.}

One measure of the complexity of branching is
the number of branching leaves between any $E, L \in {\cal E}$.
A priori there could be infinitely many in between branching
leaves producing a very complicated structure. We prove:

\vskip .1in
\noindent {\bf Theorem B} \
Let $\Phi$ be an Anosov flow in $M^3$. Let $F$ be a branching 
leaf of $\fns$ and ${\cal E}$ be the set of non separated leaves
from $F$ with the total order defined above.
Then either

(1) ${\cal E}$ is finite, hence order isomorphic to $\{1, 2, ... ,n \}$ or,

(2) ${\cal E}$ if infinite and order isomorphic to the set of integers
${\bf Z}$.

\noindent
In particular given any $E, L \in {\cal E}$, there are only finitely
many branching leaves between them.
\vskip .1in

As in the case of theorem $A$,
there is a rigidity proof of this result. However
it is quite long and complicated. Our tactic will be to first
show:

\vskip .1in
\noindent {\bf Theorem C} \
Let $\Phi$ be an Anosov flow in $M^3$ and  let $(F,L)$
be a branching pair of $\fns$. Let $g$ be a non trivial
covering translation with $g(F) = F$ 
and so that $g$ preserves transversal orientations
to $\fns, \fnu$. Then $g(L) = L$.
\vskip .1in

Using the important idea of lozenges (see definition in section $3$)
and a key result from \cite{Fe4}, theorem $B$ is an easy
consequence of theorem $C$, except that to rule out the case that
${\cal E}$ is order isomorphic to the natural numbers ${\bf N}$
we need theorem $E$ below.
Section $4$ contains a more detailed description of the set ${\cal E}$.

Theorem $C$ implies that 
$\pi(F)$ and $\pi(L)$ contain
freely homotopic closed orbits of the flow $\Phi$,
which highlights
the pervasivines of freely homotopic orbits.
This shows that the topological structure of the foliations
is intimately related to the dynamics of the flow:

\vskip .1in
\noindent {\bf Corollary D} \
Let $\Phi$ be an Anosov flow in $M^3$.
Then $\Phi$
is topologically conjugate 
to a suspension of an Anosov diffeomorphism of the
torus if and only if there are no freely homotopic
closed orbits of $\Phi$ (including non trivial free homotopies
of a closed orbit to itself).
\vskip .1in

This result does not assume that $\Phi$ is not
$\rrrr$-covered.
Another consequence of theorem $C$ is the following:

\vskip .1in
\noindent {\bf Theorem E} \
Let $\Phi$ be a non $\rrrr$-covered Anosov flow in $M^3$.
Then up to the action of covering translations, there are
finitely many branching leaves in
$\fns$.
Equivalently there are finitely many distinguished closed
orbits of $\Phi$ in $M$ so that their stable leaves
lift to branching leaves in the universal cover.
\vskip .1in

It is very important to stress here that in the above 
results we do not assume that the flow is transitive
nor is there any assumption on the manifold.
Consequently these results are the most general 
possible.
We also remark that theorems $A, B,  C$ and $E$ were 
previously proved under the assumption that $M$ has negatively
curved fundamental group and furthermore that  the flow is
quasigeodesic \cite{Fe4}. This last hypothesis 
is a very strong assumption.
The above results use only the topological
structure of the lifted foliations and have no metric 
assumption.

We also show that the structure of 
branching is strongly related to the topology of the 
ambient manifold.
We say that there is {\em infinite branching} if there are
infinitely many leaves which are not separated from each other,
otherwise we say that the branching is finite.
An easy corollary of theorem $E$ is the following:

\vskip .1in
\noindent {\bf Corollary F} \
Let $\Phi$ be an Anosov flow in $M^3$ orientable, 
atoroidal.
Then infinite branching cannot occur.
\vskip .1in

Even though the proof of corollary F is easy, it depends on
a deep result of Gabai, namely the general torus theorem
\cite{Ga} which in turn depends on the solution of the
Seifert fibered conjecture. In addition the proof uses the
characterization of Anosov flows in Seifert fibered
$3$-manifolds \cite{Gh}.
In  section $5$ we study product regions
(see definition in section $5$)  and then prove 
the following stronger result, using only the study of the topological 
structure of $\fns, \fnu$:

\vskip .1in
\noindent {\bf Theorem G} \
Let $\Phi$ be an Anosov flow in $M^3$ orientable so that there
is infinite branching in $\fns$.
Then there is infinite branching in $\fnu$.
Furthermore there is an embedded torus $T$
transverse to $\Phi$, hence $T$ is incompressible.
\vskip .1in

We remark that infinite branching does occur, for
example in the Bonatti-Langevin flow.
Furthermore we show that finite (but non trivial) branching
also occurs for a large class of Anosov flows,
for example in the flows constructed
by Franks and Williams \cite{Fr-Wi}.

Finally we apply these results to the case when $M$ has
negatively curved fundamental group.
Then $\mi$ is compactified with a sphere at infinity
$\si$. Furthermore the intrinsic geometry of a leaf
$F$ of $\fns$ or $\fnu$ is always negatively curved
in the large so there is an intrinsic ideal boundary
$\pin F$. 
In these manifolds it is fundamental to understand asymptotic
behavior of sets in $\mi$
\cite{Th1,Th2}, \cite{Mor}, \cite{Bon}.
We say that $\wwp$ has the {\em continuous extension
property} if the embedding $\varphi: F \rightarrow \mi$
extends continuously to $\varphi : 
F \cup \pin F \rightarrow \mi \cup \si$, for any leaf 
$F$ in $\fns$ or $\fnu$.
This relates the foliation to the geometry in the large of the universal
cover. This property can be defined for any Reebless
codimension $1$ foliation in  such manifolds and
it is true for fibrations \cite{Ca-Th} and many
depth one foliations \cite{Fe1}.
Recall that the  {\em limit set} of $B$ is the set of
accumulation points of $B$ in $\si$.
In this article we use the structure
of branching to analyse limit sets of leaves
when the continuous extension property holds.

\vskip .1in
\noindent {\bf Theorem H} \
Let $\Phi$ be an Anosov flow in $M^3$ with
$\pi_1(M)$ negatively curved.
Suppose that $\wwp$ has the continuous extension
property. If $\Phi$ is not $\rrrr$-covered
then the limit set of any leaf
$F$ of $\fns$ or $\fnu$ is a Sierpinski curve,
that is, the complement of a countable, dense
union of open disks in $\si$.
\vskip .1in

If $\Phi$ is $\rrrr$-covered then the limit
set of any $F \in \fns \cup \fnu$ is $\si$,
regardless of whether the continuous
extension property holds or not \cite{Fe2}.
There are many $\rrrr$-covered examples
\cite{Fe3}. 

Bonatti and Langevin's example of a transitive, 
non $\rrrr$-covered Anosov flow in dimension
$3$, was generalized by Brunella \cite{Br}
who produced many examples by Dehn surgery
on geodesic flows. The tool used to show
that these flows are not $\rrrr$-covered was
the existence of a transverse torus; hence
all such examples were not in hyperbolic
$3$-manifolds. The main open conjecture in
this theory was whether
$M$ being hyperbolic would imply that
the flow $\Phi$ is $\rrrr$-covered.
In this article we
answer this conjecture in the negative:

\vskip .1in
\noindent {\bf Theorem I} \
There is a large class of transitive, 
non $\rrrr$-covered Anosov flows
where the underlying $3$-manifold is hyperbolic.
This includes all Anosov flows in 
non orientable, hyperbolic $3$-manifolds.
\vskip .1in

In a forthcoming paper \cite{Fe8} we use the results of this 
article to study incompressible tori in $3$-manifolds supporting
Anosov flows. It is of great interest to find, in the isotopy
class of the torus, the best position with respect to the flow
\cite{Ba3,Ba4}. We prove:

\vskip .1in
\noindent {\bf Theorem} \ (\cite{Fe8})
Let $\Phi$ be an Anosov flow in $M^3$ and let $T$ an incompressible
torus in $M$. Suppose that no loop in $T$ is freely homotopic
to a closed orbit of $\Phi$. Then $\Phi$ is topologically 
conjugate to a suspension Anosov flow.
Furthermore $T$ is isotopic to a torus transverse to $\Phi$.
\vskip .1in

The article is organized as follows: in the next section
we develop background material. In section $3$ we prove
theorem $A$ and in the following section we prove theorem $C$
and immediately derive theorems $B$ and $E$ and 
corollaries  $D$ and $F$. Section $5$ studies product
regions, which is then applied to a more detailed
analysis of infinite branching 
and the construction of a transverse torus in section $6$. In
the final section we  study the continuous extension property.

We thank Bill Thurston for 
encouragement and many helpful conversations relating to this work.
We also thank Thierry Barbot for useful suggestions to a first 
version of this article.

\section{Background}

Let $\Phi _t : M \rightarrow M$ be a nonsingular flow in a closed 
Riemannian manifold $M$.
The flow $\Phi$ is {\em Anosov} if there is a continuous decomposition
of the tangent bundle $TM$ as a Whitney sum  
$TM = E^0 \oplus  E^s \oplus  E^u$ of $D \Phi _t$ invariant subbundles
and there are constants $\mu_0 \geq 1$, $\mu_1 > 0$ so that:

(i) $E^0$ is one dimensional and tangent to the flow,

(ii) $|| D \Phi _t (v) || \leq  \mu_0 e^{- \mu_1 t} || v ||$  
\ for any $v \in E^s, \ t \geq 0$,

(iii) $|| D \Phi _{-t} (v) || \leq \mu_0 e^{- \mu_1 t} || v ||$  
\ for any $v \in E^u, \ t \geq 0$.

In this article we restrict to $M$ of dimension $3$. 
Then $E^s, E^u$ are one dimensional
and integrate to one dimensional foliations $\fss, \fuu$ 
called the strong stable
and strong
unstable foliations of the flow. Furthermore, the bundles 
$E^0 \oplus  E^s$ and $E^0 \oplus  E^u$ are also integrable 
\cite{An} producing
$2$-dimensional foliations $\fs, \fu$ 
which are the stable and unstable foliations
of the flow. 

The flow is said to be {\em orientable} when both
$\fs, \fu$ are transversely orientable.
We remark that there is always a regular cover
of order $\leq 4$ where the lifted $\fs$ and $\fu$
are transversely orientable. Whenever 
possible we will lift to such a cover.

The leaves of $\fs, \fu$ are either
topological planes, annuli or M\"{o}ebius bands. 
The last two correspond exactly to leaves containing closed 
orbits of $\Phi$.
There is at most one closed orbit in a leaf of $\fs$,
in which case all other orbits are forward 
asymptotic to it. Similarly for $\fu$.


The foliation $\fs$ is Reebless, so Novikov's
theorem \cite{No} implies that
given any closed orbit $\gamma$ of $\Phi$,
$\gamma^n$ is not null homotopic for any $n \not 
= 0$.

Let 
$\pi: \widetilde M \rightarrow M$
be the universal covering
space of $M$. 
This notation will be fixed throughout the article.
The Anosov foliations $\fs, \fu$ lift to foliations
$\fns, \fnu$ in $\widetilde M$.
The leaves of $\fns, \fnu$ are topological planes, so
$\widetilde M$ is homeomorphic to $\rrrr ^3$ \cite{Pa}.
Therefore $M$ is {\em irreducible}
that is every embedded sphere in $M$ bounds a $3$-ball.
The induced flow in $\widetilde M$ is denoted by $\wwp$. 

Let $\oo$ be the orbit space of $\wwp$ obtained 
by collapsing flow lines to points
and let 
$\Theta: \widetilde M \rightarrow \oo$
\noindent
be the projection map.
A fundamental property which will be repeatedly used here is
that $\oo$ is Hausdorff and hence homeomorphic to $\rrrr^2$
\cite{Fe3}.
This is a significant simplification since now much of
the analysis can be done in dimension $2$ instead of dimension $3$.
We stress that $\oo$ is only a topological
object. There is no natural metric in $\oo$ since the flow
direction contracts and expands distances in $\mi$.
The foliations $\fns, \fnu$ induce two transverse $1$-dimensional
foliations in $\oo$, which will also be denoted by $\fns, \fnu$.
By an abuse of notation we  will many times identify 
sets in $\mi$ or orbits of $\wwp$ 
to their respective images in ${\cal O}$.

The fundamental group $\pi_1(M)$ is isomorphic
to the set of covering translations of
$\mi$. We will usually assume one such
identification is fixed.
Given a covering translation $g$, we will
also denote by $g$ its action on $\hs, \hu$.


Let $W^s(x)$ be the leaf of $\fs$ containing $x$ and similarly
define $W^u(x)$, $W^{ss}(x), W^{uu}(x), \widetilde W ^s(x),
\widetilde W ^u(x), \widetilde W ^{ss}(x)$ and $\widetilde W ^{uu}(x)$.
In the same way if $\alpha$ is an orbit of $\Phi$ we define
$W^s(\alpha)$, etc..
General references for Anosov flows are 
\cite{An}, \cite{An-Si}, \cite{Bow}, \cite{Sh} 
and \cite{Sm}.

An {\em incompressible} surface ($\not = {\bf S}^2$)
is an embedded surface
in $M^3$ which is injective in the fundamental group level.
A manifold is {\em toroidal} if it contains an
incompressible torus and {\em atoroidal} otherwise.

\section{Periodic branching leaves}

The following definitions will be useful.
If $L$ is a 
leaf of $\fns$ or $\fnu$,
then a {\em half leaf} of $L$ is a connected component
$A$ of $L - \gamma$, where $\gamma$ is any
full orbit in $L$. 
The closed half leaf is $\overline A = A \cup \gamma$ and 
its boundary is $\partial A = \gamma$.
If $L$ is a 
leaf of  $\fns$ or $\fnu$ then a {\em flow band} $B$
defined by orbits $\alpha \not = \beta$ in $L$
is the connected component of $L - \{ \alpha, \beta \}$
which is not a half leaf of $L$.
The closed flow band
associated to it is $\overline B =
B \cup \{ \alpha, \beta \}$ and its boundary is
$\partial B = \{ \alpha, \beta \}$.

Since $\mi$ is simply connected, $\fns$ and $\fnu$
are transversely orientable. 
Choose one such orientation, 
assumed to agree with
the lifts of the  transversal orientations to $\fs, \fu$ if
any of these is transversely oriented. 
Notice that in general, covering translations may not
preserve transversal orientations.

For $p \in \widetilde M$,
let $\widetilde W^s_+(p)$ be the half
leaf of $\widetilde W^s(p)$ defined by the
orbit $\wwp _{\rrrr} (p)$
and the positive transversal orientation to $\fnu$
at $p$. This is also called 
a {\em positive half leaf} of $\ws(p)$.
Similarly define $\widetilde W^s_-(p)$, 
a {\em negative half leaf} and also define
$\widetilde W^u_+(p)$ and 
$\widetilde W^u_-(p)$.

A fundamental property for us is that
any leaf $L$ in $\fns$ or $\fnu$ separates $\mi$.
This is a consequence of $\mi$ being simply connected.
The {\em front} of $L$ is the component
of $\widetilde M - L$ defined by the
positive  transversal orientation to $L$.
Similarly define the {\em back} of $L$.
For $p \in \mi$ let $\widetilde W^{ss}_+(p) =
\widetilde W^s_+(p) \cap \widetilde W^{ss}(p)$.
In the same way define 
$\widetilde W^{ss}_-(p)$,
$\widetilde W^{uu}_+(p)$ and
$\widetilde W^{uu}_-(p)$.

If $F \in \fns$ and $G \in \fnu$ 
then $F$ and $G$ intersect in at most one
orbit, since two intersections would 
force a tangency of $\fns$ and $\fnu$.
This is easiest seen in $\oo$, as $\fns$ and
$\fnu$ are then $1$-dimensional foliations
of the plane.

We say that leaves $F, L \in \fns$ and $G, H \in \fnu$
form a {\em rectangle}  if $F$ intersects both 
$G$ and $H$ and so does $L$, see fig. \ref{bands} a.
We also say that $E$ intersects
$G$ {\em between} $F$ and $L$ if $E \cap G$ is contained
in the flow band in $G$ defined by $G \cap F$ and
$G \cap L$.
Then it is easy to prove \cite{Fe5} that if $E \in \fns$
intersects $G$ between $F$ and $L$ then $E$ also
intersects $H$ between $F$ and $L$.
This means that there is  a product structure of
$\fns$ and $\fnu$ in the region bounded
by $F, L, G$ and $H$.

The following two definitions will be essential for all resuts in
this article.

\begin{define}{}{}
Given $p \in \mi$ (or $p \in \oo$), let 

$$\ooo^u_+ (p)  = \{ F \in \fns  \ | \ F \cap 
\widetilde W^u_+(p) \not = \emptyset \},$$

\noindent
an open subset of $\hs$.
Notice that the leaf $\ws(p) \not \in \ou(p)$.
Similarly define $\ooo^u_-(p), \ooo^s_+(p)$ and $\ooo^s_-(p)$.
%
%
\end{define}

\blankfig{bands}{1.8}{a. Rectangles,
b. Perfect fits in the universal cover.}

\begin{define}{}{}
Two leaves $F, G$, $F \in \fns$ and $G \in \fnu$, form
a perfect fit if $F \cap G = \emptyset$ and there
are half leaves $F_1$ of $F$ and $G_1$ of $G$ 
and also flow bands $L_1 \subset L \in \fns$ and
$H_1 \subset H \in \fnu$, (see figure \ref{bands} b)  so that:

$$ \overline L_1 \cap \overline G_1 = \partial L_1 \cap \partial G_1,
\ \ \overline L_1 \cap \overline H_1 = \partial L_1 \cap \partial H_1,
\ \ \overline H_1 \cap \overline F_1 = 
\partial H_1 \cap \partial F_1,$$ 

$$\forall \ S \in \fnu, \ \ \ \   
S \cap L_1  \not = \emptyset \ \Leftrightarrow
S \cap F_1  \not = \emptyset \ \ {\rm and}$$

$$\forall \ E \in \fns, \ \ \ \ E \cap G_1  
\not = \emptyset \ \Leftrightarrow
E \cap H_1  \not = \emptyset.$$

\end{define}

\noindent
Notice that the flow bands $L_1, H_1$ 
(or the leaves $L, H$) are not uniquely 
determined given the perfect fit $(F,G)$.
We will also say that $F$ and $G$ are {\em asymptotic} in the sense
that if we consider stable leaves near $F$ and on the
side containing $G$ they will intersect $G$ and vice versa.
Perfect fits produce ``ideal" rectangles, in the sense that
even though $F$ and $G$ do not intersect, there is
a product structure (of $\fns$ and $\fnu$) in the 
interior of the region bounded by $F, L, G$ and $H$.

It is easy to show \cite{Fe5} that there is at most one
leaf $G \in \fnu$ making a perfect fit with a
given half leaf of $F \in \fns$ and in 
a given side of $F$.
Therefore a perfect fit
is a detectable property in $M$.
This means that if $(L,G)$ forms a perfect fit and $g$ is
any orientation preserving covering translation with $g(L) = L$,
then $g(G) = G$. The last assertion follows from
uniqueness of perfect fits and the fact
that, as $g$ acts by homeomorphisms in the leaf spaces,
it takes perfect fits to  perfect fits.

%

If $p, q$ are in the same strong stable leaf
let $[p,q]_s$ denote the closed segment in that leaf
from $p$ to $q$ and let $(p,q)_s$ be the corresponding
open segment. Similarly define $[p,q]_u$ and $(p,q)_u$.

We say that $\ooo^s_+(p)$ and $\ooo^s_+(q)$ are
comparable and  will denote this  by 
$\ooo^s_+(p) \sim  \ooo^s_+(q)$, if one of them is
contained in the other. Then we write
$\ooo^s_+(p) < \ooo^s_+(q)$ if the former is
strictly contained in the latter. 
Similarly define $\leq, >$ and $\geq$.
The symbol $\not \sim$ means not comparable.

We also say that an orbit $\gamma$ of $\wwp$ is 
{\em periodic} if it is left invariant by a non
trivial covering translation.

\begin{theorem}{}{}
Let $\Phi$ be an Anosov flow in $M^3$ 
and let $F$ be a branching leaf of $\fns$. Then there 
is a non trivial covering translation $g$ with $g(F) = F$,
that is, $F$ is periodic.
\label{peri}
\end{theorem}

\begin{proof}{}
By taking a finite cover if necessary,
we may assume that $\Phi$ is orientable. 
Let $L \in \fns$, $L \not = F$, so that
$F,L$ form a branching pair of $\fns$.
Assume without loss of generality
that $F$ and $L$ are not separated on
their negative sides, that is they
are associated to branching of $\fns$ in
the positive direction (positive branching).

Let $w_0 \in F$, $w' \in L$. Since $F$ and $L$
are not separated in their negative sides
there are $y_0 \in \widetilde W^{uu}_-(w_0)$
($y_0$ sufficiently near $w_0$)
and $x_0 \in \wu(w') \cap \widetilde W^{ss}(y_0)$
so that if $r_0 = \widetilde W^{uu}_+(x_0) \cap L$,
then for any $E \in \fns$,

$$E \cap (y_0,w_0)_u \not = \emptyset \ \ \Leftrightarrow
\ \ E \cap (x_0,r_0)_u \not = \emptyset. \ \ \ \ (*)$$

\noindent 
This fact,
which follows from the separation property of
leaves of $\fns$,
will often be implicitly used.

By switching $F$ and $L$ if necessary we may
assume that $\wu(x_0)$ is in the front of $\wu(y_0)$.
Our first goal will be to find unique leaves associated to
the branching which
form perfect fits with $F$ and $L$.

As there are $z \in [y_0,x_0]_s$ with  $\wu(z) \cap F = \emptyset$
(for instance $z = x_0$),
let $p_0$ be the closest point to $y_0$ in $[y_0,x_0]_s$ so
that $\wu(p_0) \cap F = \emptyset$.

\begin{lemma}{}{}
The  leaves
$F$ and $\wu(p_0)$ form a perfect fit.
\label{fit}
\end{lemma}

\begin{proof}{}
For candidates of flow bands 
let $A = \wwp _{\rrrr}((y_0,w_0)_u)$ and
$B = \wwp _{\rrrr}((y_0,p_0)_s)$.
Then $\overline A \cap \overline B =
\wwp_{\rrrr}(y_0)$, 
$\overline A \cap F = \wwp_{\rrrr}(w_0)$
and $\overline B \cap \wu(p_0) = \wwp _{\rrrr}(p_0)$.

Let $E \in \fns$ with $E \cap A \not = \emptyset$.
Then $E \cap \wu(x_0) \not = \emptyset$.
Since $\wu(p_0)$ separates $M$ it follows that
$E \cap \wu(p_0) \not = \emptyset$.
As $E$ is in front of $\ws(y_0)$ then 
$E \cap \widetilde W^u_+(p_0) \not = \emptyset$.

\blankfig{split}{1.7}{Branching in $\fns$.}

Conversely let $E \in \fns$ with 
$E \cap \widetilde W^u_+(p_0) \not = \emptyset$.
Suppose that  $E \cap A = \emptyset$. Since
$\wu(p_0) \cap F = \emptyset$, then the front
of $E$ is disjoint from the front of $F$. For any
$z \in \widetilde W^{ss}(p_0)$ near enough $p_0$,
$\wu(z) \cap E \not = \emptyset$. 
As $E$ is in the back of $F$, it follows that 
$\wu(z) \cap 
F = \emptyset$. 
This contradicts the choice of $p_0$.
We conclude that $E \cap A \not = \emptyset \Leftrightarrow
E \cap \widetilde W^u_+(p_0) \not = \emptyset$.

\vskip .1in
Let now $R \in \fnu$ with $R \cap B \not = \emptyset$.
If $R \cap F = \emptyset$, then $z = R \cap [y_0,p_0]_s$
is closer to $y_0$ (in $\widetilde W^{ss}(y_0)$)
than $p_0$, contradiction. Hence $R \cap F \not
= \emptyset$,
in particular $R \cap \widetilde W^s_+(w_0) \not = \emptyset$.

Conversely suppose that $R \cap \widetilde W^s_+(w_0) 
\not = \emptyset$. 
Let $F^* \in \fns$ be close enough to $F$
so that $F^* \cap R \not = \emptyset$,
$F^* \cap \wu(y_0) \not = \emptyset$ and
$F^* \cap \wu(x_0) \not = \emptyset$.
Then $\wu(y_0), \wu(x_0), \ws(y_0)$ and $F^*$ form 
a rectangle. 
Since and $R \cap F^* \not = \emptyset$ is 
between $F^* \cap \wu(y_0)$
and $F^* \cap \wu(x_0)$ 
then $R \cap \ws(y_0) \not = \emptyset$.
As $R$ is in front of $\wu(y_0)$
then $R \cap \widetilde W^s_+(y_0) \not = \emptyset$.
Since $R \cap F \not = \emptyset$
then $R$ is in the back of $\wu(p_0)$.
Therefore $R \cap B \not = \emptyset$.
This finishes the proof of the lemma.
\end{proof}

\noindent
\underline{Continuation of the proof of theorem \ref{peri}}

In the same way there is a unique $q_0 \in [y_0,x_0]_s$
with $\wu(q_0)$ and $L$ forming a perfect fit.
By uniqueness of perfect fits,
the leaves $\wu(p_0), \wu(q_0)$ depend only on
$F$ and $L$. 
If follows from $(*)$ and lemma
\ref{fit}, that given $E \in \fns$, 
$E \cap \widetilde W^u_+(p_0) \not = \emptyset
\Leftrightarrow E \cap \widetilde W^u_+(q_0) \not = \emptyset$.
Equivalently $\ou(p_0) = \ou(q_0)$.

\vskip .1in
\noindent \underline {Case 1} \
$p_0 = q_0$.

Let $G = \wu(p_0) = \wu(q_0)$. If $G$ is periodic there is
$g \not = id$ with $g(G) = G$. By uniqueness of perfect
fits and preserving of transversal orientations it
follows that $g(F) = F$
and we are done. So we may assume that $G$ is not periodic.

Let $c_0 = \pi(p_0)$. Since $G$ is not periodic,
$\Phi_{\rrrr}(c_0)$ is not a closed orbit, nor is
it backwards asymptotic to a closed orbit. Let $c$ be a negative
limit point of $\Phi_{\rrrr}(c_0)$
and let $c_i = \Phi_{t_i}(c_0)$, $t_i \rightarrow -\infty$,
with $c_i \rightarrow c$.
If $c_i$ and $c_j$ are in the same local unstable
leaf near $c$, then there is a closed path in
$W^u(c_i)$ consisting of the flow segment from $c_i$
to $c_j$ and then a small strong unstable segment
from $c_j$ to $c_i$ in the local unstable leaf
through $c_j$.
This path is not null homotopic in $W^u(c_i)$,
hence $W^u(c_i)$ contains a closed orbit,
contradiction to our assumption. This is the key fact used
in the proof of the theorem and it will imply that non periodic
leaves in the universal cover are not rigid .

Lift $c_i$ to $p_i \in \mi$ with $p_i
\rightarrow p$ and $\pi(p) = c$.
Then $p_i = g_i(\wwp_{t_i}(p_0))$, where $g_i$ are
covering translations.
By the above argument $\ws(p_i) \not = \ws(p_k)$
for any $i \not = k$. This is the non rigidity we
are looking for.

\blankfig{pert}{2.4}{Rigidity of branching leaves: the adjacent case}

Let $F_i = g_i(F), L_i = g_i(L), A_i = g_i(A)$,
$B_i = g^i(B)$
and $G_i = g_i(G)$.
Let $y_i = g_i(\wwp_{t_i}(y_0))$
and let $x_i = g_i(\wwp_{t_i}(x_0))$.
Up to subsequence
assume that all $p_i$ and $p$ are near enough, in a product 
neighborhood of $\fnu$ of diameter $<< 1$. 
Assume also that for all $i$,

$$l(\wwp_{t_i}([y_0,p_0]_s)) > 1  \ \ {\rm and} \ \
l(\wwp_{t_i}([p_0,x_0]_s)) > 1. \ \ \ \ (**)$$   

Choose $i, k$ so that $p_i$ is in the back of $\widetilde W^u(p_k)$,
see fig. \ref{pert}.
Since $d(p_i, p_k) << 1$ it follows that  
$\widetilde W^s_-(p_k)
\cap \wu(p_i) \not = \emptyset$
and $\widetilde W^s_+(p_i)
\cap \wu(p_k) \not = \emptyset$. 
By $(**)$, this implies that
$y_k$ is in the back of $\wu(p_i)$
and $x_i$ is in the front of $\wu(p_k)$, see fig. \ref{pert}.
Hence $\wu(y_k)$ is in the back of $\wu(p_i)$.
Then $\wu(p_i) \cap B_k \not = \emptyset$,
hence $\wu(p_i) \cap F_k \not = \emptyset$.
As $L_i$ makes a perfect fit with $\wu(p_i)$,
this implies that $L_i$ is in front of $F_k$,
hence $L_i$ is in the back of $\wu(p_k)$.

On the other hand, $L_i \cap \wu(x_i) \not = \emptyset$.
Since $\wu(x_i)$ is in front of $\wu(p_k)$
then $\wu(p_k) \cap \wwp _{\rrrr}([p_i,x_i]_s) \not = \emptyset$.
As $L_i$ and $\wu(p_i)$ form  a perfect fit, this implies that
$\wu(p_k) \cap L_i \not = \emptyset$. This contradicts 
the previous paragraph.

This shows that if $p_0 =  q_0$,
then $G$ is periodic, left invariant
by $g$, hence $F$ and $L$ are periodic
and both left invariant under $g$.

\vskip .1in
\noindent
{\bf Remarks:} 
(1) If we apply the argument above
when $G$ is periodic, we get $\ws(p_i) = \ws(p_k)$
for all $i, k$. There is no small perturbation
of the local picture, which is then rigid.
This will imply that the whole set of non  separated leaves
from $F$ is very  rigid.

(2) It is tempting to try the following ``intuitive"
approach to  the above proof: as $\pi(\widetilde W^u_+(p_0))$
is not compact in $M$, there are always translates $S_1$ and
$S_2$ of $\wu(p_0)$ and points $u_i \in S_i$ arbitrarily
near  each other.  However there is no control of the rest
of the picture. For instance we do not know a priori
what happens to the respective stable lengths.
This is the reason why we fixed an orbit $\Phi_{\rrrr}(\pi(p_0))$
and flowed {\underline {backwards}}  in order to insure
that stable lengths are as big as we want.
\vskip .1in

\blankfig{diff}{1.9}{Rigidity of branching: the separated case.}

\vskip .1in
\noindent \underline {Case 2} \
$p_0 \not = q_0$.


We use the same notation as in case $1$.
As $q_0  \not = p_0$, let 
$q_i = g_i(\wwp_{t_i}(q_0))$.
Choose $i, k$ with $p_i$ in the back of
$\wu(p_k)$.
As in case $1$, 
$\widetilde W^u_+(p_i) \cap F_k \not = \emptyset$. 
There is no a priori contradiction because
now $L_i$ does not form a perfect
fit with $\wu(p_i)$, and in fact $L_i$ is probably in 
the front of $\wu(p_k)$.
Let 

$$e_1 = \wu(p_k) \cap \widetilde W^{ss}_+(p_i), \ \ 
e_2 = \wu(p_i) \cap \widetilde W^{ss}_-(p_k).$$

\noindent
Then $\ou(p_k) < \ou(e_2)$ and by the local
product structure of $\fns, \fnu$ near $p$, it
follows that 
$\ooo ^u_+ (p_i) > \ooo ^u_+ (e_1)$, see fig. \ref{diff}.
Choose $E \in 
\ooo ^u_+ (p_i) - \ooo ^u_+ (e_1)$.
By the above considerations it is clear that
$E \cap \wu(p_k) = \emptyset$.
But 

$$\ou(q_i) = \ou(g_i(\wwp_{t_i}(q_0))) = 
g_i(\ou(\wwp_{t_i}(q_0))) = g_i(\ou(\wwp_{t_i}(p_0)))
= \ou(p_i),$$ 

\noindent
hence $E \in \ou(q_i)$.
As a result $E \cap \widetilde W^u_+(q_i)
\not = \emptyset$. 
But now $\wu(q_i)$ is in the front of $\wu(p_k)$. 
Since $\wu(p_k)$ separates 
$\mi$, then $E \cap \wu(p_k) \not = \emptyset$,
contradiction. As before we conclude that $G$ is periodic, left
invariant by $g \not = id$, so
$F$ is also left invariant by $g$. 
\end{proof}

\noindent {\bf Caution:}
The same argument shows that $L$ and
$\wu(q_0)$ are also periodic.
We do not know at this point that the 
same covering translation
leaves invariant {\underline {both}} $F$ and $L$.
This is a much stronger fact.

\section{Branching structure}

In this section we show that if
$F$ and $L$ are not separated,  then not only 
are they periodic, but there is a common covering translation
leaving both of them invariant.
As a result, branching forces 
a non trivial free homotopy between closed
orbits of $\Phi$ in $M$ and this gives the topological
characterization of suspensions.
Furthermore we will show that $F$ and $L$ are 
connected by a finite
sequence of lozenges, as defined below.
This completely determines the structure of the set
of non separated leaves from $F$.
As a consequence we show there are only
finitely many branching leaves up to
covering translations.
This in turn implies that if there is infinite branching then there
is an incompressible torus in $M$.

\begin{define}{}{}
Lozenges - 
Let $p, q \in \widetilde M$, $p \not \in \ws(q)$,
$p \not \in \wu(q)$. Let $H_p$ be the half leaf
of $\wu(p)$ defined by   $\wwp _{\rrrr} (p)$ and
contained in the same side of $\ws (p)$ as $q$.
Let $L_p$ be the similarly defined half leaf of $\ws(p)$ and
in the same fashion define  $H_q, L_q$.  Then $p, q$ form 
a lozenge, fig. \ref{loz}, a if
$H_p, L_q$ and $H_q, L_p$ respectively form
perfect fits.
\end{define}

\blankfig{loz}{2.1}{a. A lozenge, b. A chain of adjacent lozenges.}

We say that $p, q$ (or $\wwp_{\rrrr}(p),
\wwp_{\rrrr}(q)$) are corners of the lozenge.
If the lozenge with corner $p$ is contained
in the back of $\ws(p)$ then $p$ is a corner
of type $(+,*)$, otherwise it is of type $(-,*)$.
Similarly using $\wu(p)$ define types $(*,+), (*,-)$.
The sides of the lozenge are $H_p, L_p, H_q$ and $L_q$.
Since given any four leaves there is at  most one lozenge
defined by them we will also say the full leaves
are the sides of the lozenge.
Notice that if $p$ is a corner of type $(-,-)$ 
then $\ou(p) = \ooo^u_-(q)$, $\ooo^s_+(p) = \ooo^s_-(q)$ 
 and similarly for the  other cases.

Two lozenges are {\em adjacent} if they share a corner and
there is a stable or unstable leaf
intersecting both of them, see fig. \ref{loz} b.
A chain of lozenges is a collection $\{ {\cal B}_i \}, 
1 \leq i \leq n$, where $n \in {\bf N}$, so that 
${\cal B}_i$ and ${\cal B}_{i+1}$ share
a corner. Consecutive lozenges may be adjacent
or not.

The following theorem will be essential for the results
in this section:

\begin{theorem}{Fe4}{}
Let $\Phi$ be an Anosov flow in $M^3$.
Suppose that $F_i, i = 0,1$ are leaves of $\fns$ for
which there is a non trivial covering translation $g$
with $g(F_i) = F_i, i = 0,1$.
Let $\alpha_i, i = 0,1$ be the periodic orbits of $\wwp$
in $F_i$ so that $g(\alpha_i) = \alpha_i$.
Then $\alpha_0$ and $\alpha_1$ are connected
by a finite chain of lozenges 
$\{ {\cal B}_i \}, 1 \leq i \leq n$ and $g$
leaves invariant each lozenge 
${\cal B}_i$ as well as their corners.
\label{chain}
\end{theorem}

Furthermore there is a unique chain that is
minimal, in the sense that any other chain from
$\alpha_0$ to $\alpha_1$ contains this
chain \cite{Fe7}.
Given any chain ${\cal B} = \{ {\cal B}_i \}, 1 \leq i \leq n$
from $\alpha_0$ to $\alpha_1$, 
let $\gamma_0 = \alpha_0$ and inductively
define $\gamma_i, i > 0$ to be the remaining corner
of ${\cal B}_i$.
The minimal chain from $\alpha_0$ to $\alpha_1$ is defined by:
${\cal B}_{i+1}$ is on the same
side of $\ws(\gamma_i)$
and $\wu(\gamma_i)$ that $\alpha_1$ is.

A closed orbit of $\Phi$ traversed once 
is called an {\em indivisible} closed orbit.

The following result will be often used 
in this article:

\begin{theorem}{Fe7}{}
Let $\Phi$ be an orientable Anosov flow in $M^3$.
If $\gamma$ is an indivisible closed orbit
of $\Phi$, then $\gamma$ represents
an indivisible element in $\pi_1(M)$.
Equivalently if $g^n(F) = F$, where
$F \in \fns \cup \fnu$, $g$ is a covering
translation and $n \not = 0$, then $g(F) = F$. 
\label{indi}
\end{theorem}

There is a related result if $\Phi$ is not
assumed to be orientable.

The {\em stabilizer} ${\cal T(F)}$ of a leaf $F$ of $\fns$ (or $\fnu$)
is the subgroup of $\pi_1(M)$ of those $g$ with $g(F) = F$.
If $\pi(F)$ does not contain a periodic orbit, then
${\cal T}(F)$ is trivial. 
Otherwise let $\gamma$ be the indivisible closed orbit in
$\pi(F)$. Then ${\cal T}(F)$ is infinite cyclic and
it has a generator conjugate to $[\gamma]$ in $\pi_1(M)$.

The main technical result in this section is the following:

\begin{theorem}{}
Let $\Phi$ be an Anosov flow in $M^3$. Suppose that $F, L$
form a branching pair of $\fns$. 
Let $g$ be a non trivial covering translation
with $g(F) = F$, so that $g$ preserves
transversal orientations to $\fns, \fnu$. Then $g(L) = L$.
Similarly for $\fnu$.
\label{theb}
\end{theorem}

\begin{proof}{}
Up to a finite cover  assume that $\Phi$ is orientable.
Since $g$ preserves transversal orientations,
then $g$ is still a covering translation of the finite cover.
Without loss of generality suppose that  $F$ and $L$
are not separated on  their negatives sides,
corresponding to positive branching.
Finally we may assume that $g$ generates ${\cal T}(F)$.

As in theorem \ref{peri} there are unique 
leaves $G, H \in \fnu$ making perfect
fits with $F$ and $L$ respectively and so that:
$G$ separates  $F$ from $L$  and  so does $H$.
Let $p \in G$ so that $\widetilde W^{ss}(p)$
intersects $H$ and let $q = \widetilde W^{ss}(p)
\cap H$.
Recall from the proof of theorem \ref{peri}, that
$\ou(p) = \ou(q)$.

Since $g$ preserves transversal orientations then 
$g(G) = G$. Our goal is to show that $g(L) = L$.
Suppose then that $g(L) \not = L$, hence 
by the same argument $g(H) \not = H$. 
Let $\gamma \subset G$ be the periodic orbit of $\wwp$ 
in $G$, so 
$g(\gamma) = \gamma$. 

\vskip .1in
\noindent \underline{Claim 1}  \ -
There is $R \in \fnu$
in the back of $L$ making a perfect fit  with a positive
half leaf of $L$, hence $R$ is in the front of $H$.

We may assume that $p \in \widetilde W^u_+(\gamma)$.
Let $E = \ws(p)$. 
By taking $g^{-1}$ if necessary 
assume that $g(E)$ is in front of $E$. 
Hence $g(E) \in \ooo^u_+(p)$, therefore $g(E) \in \ooo^u_+(q)$.
Then $H \cap g(E) \not = \emptyset$. 
There are $2$ cases:

\vskip .1in
(1) $g(H)$ is in front of $H$, see fig. \ref{iter}.

Let $e' = \widetilde W^{ss}(g(p)) \cap H$.
Since $g(p) \in \widetilde W^u_+(p)$
then $\ou(g(p)) = \ou (e')$.
But also $\ou(g(p)) = \ou (g(q))$,
so $\ou(g(q)) = \ou(e')$, where $g(q) \in g(H)$
and $e' \in H$.
Since $L$ makes a perfect fit with $H$ and
$g(L)$ makes a perfect fit with $g(H)$ 
this shows that $g(L)$ is not separated from $L$.

\blankfig{iter}{2.0}{Iterating non invariant leaves.}

As in the proof of theorem \ref{peri}, there is a unique
$e_0 \in [e',g(q)]_s$ with $\wu(e_0)$ making 
a perfect fit with $L$ and $L$ in the back of 
$\wu(e_0)$. In this case let $R = \wu(e_0)$.

\vskip .1in
(2) Suppose now that $g(H)$ is in the back of $H$.

Notice that $E, g(E), H$ and $G$ form a rectangle.
Since $g(H) \cap g(E) \not = \emptyset$ and
$g(H)$ is between $G$ and $H$ it follows
that $g(H) \cap E \not = \emptyset$
and $g(H) \cap E$ is an orbit in $E$ between
$E \cap G$ and $E \cap H$.

In this case let $c = g(H) \cap \widetilde W^{ss}(p)$.
Then $c \in  (p,q)_s$. 
Since $\ou(p) = \ou (q)$ and 
$L \not \in \ou(q)$, then
$g(L) \not \in \ou(q)$, 
so $g(L) \cap H = \emptyset$. Hence $g(L)$ is in the back of $H$.
As in case $(1)$, it follows that $L$ and $g(L)$ form
a branching pair.
Let $c_2 \in (c,q)_s$ with
$\wu(c_2)$ making a perfect fit
with $g(L)$ and with $g(L)$ in the back of $\wu(c_2)$.
Then $R = g(\wu(c_2))$ makes  a perfect fit with $L$
and $L$ is in the back of $\wu(c_2)$.
This finishes the proof of claim $1$.

\vskip .1in
By theorem \ref{peri}, $L$ is periodic and let 
$\alpha^*$ be the indivisible periodic orbit in $L$.
Let $h$ a generator of ${\cal T}(H)$.
Since $\Phi$ is orientable, $h(H) = H$, $h(R) = R$.
Let $\alpha$ be the periodic orbit in $H$.
Therefore
$L$ and $H$ are 2 of the sides of a
lozenge ${\cal N}_1$  with
other sides in $\ws(\alpha)$ and $\wu(\alpha^*)$,
that is $\alpha$ and $\alpha^*$ are the 
corners of the lozenge.
In the same way $L$ and $R$ are the
2 sides of a lozenge ${\cal N}_2$. 
The lozenges are adjacent and intersect the stable leaf $E$.
Let ${\cal N} = {\cal N}_1 \cup {\cal N}_2$.

We now show that $F$ also makes a perfect fit with $U \in \fnu$,
$U \not = G$ and $F$ in the front of $U$, hence $G$ is in the front
of $U$, see fig. \ref{bel}. If $h(G) = G$ then since 
$g$ generates ${\cal T}(G)$,
it follows that  $h = g^n$ for
some $n \in {\bf Z}$. Hence $g^n(H) = H$.
Theorem \ref{indi} then implies that
$g(H) = H$ contrary to assumption.
It follows that $h(G) \not = G$.
Using claim $1$ with the roles of $F, L$ exchanged, we produce 
the required $U \in \fnu$.
Furthermore there are two adjacent lozenges 
${\cal D}_1$ and ${\cal D}_2$  with (some)
sides in $U, F, G$. 
Let ${\cal D}$ be their union. Both lozenges intersect
a stable leaf which we may assume is $E$.

\blankfig{bel}{1.4}{Double lozenges.}

\vskip .1in
From now on the proof goes roughly as follows:
We will show that $\ws(\gamma)$ intersects 
of $\widetilde W^u_+(\alpha)$ and similarly
that $\ws(\alpha)$ intersects $\widetilde W^u_+(\gamma)$
producing a contradiction.

By taking $g^{-1}$ if necessary
suppose that $g(H)$ is in the back of $H$.
Let $H_i = g^i(H)$. Then as in case $(2)$ 
of the claim, $H_{i+1}$ is in the back
of $H_i$, and for all $i \geq 0$,  $H_i \cap E \not = \emptyset$.
Furthermore $H_i$ is always in front of $G$. This implies that
$H_i \rightarrow
S$ with $S \cap E \not = \emptyset$  (and maybe $H_i$ also converges to
other leaves  of $\fnu$). 

Let ${\cal A}_i$ be the front of $H_i$ and let ${\cal A} =
\cup _{i \in {\bf N}} {\cal A}_i$.
Then $g({\cal A}_i) = {\cal A}_{i+1}$ 
so $g({\cal A}) = {\cal A}$ and consequently
$g(\partial {\cal A})  = \partial {\cal A}$.
Since $S \not \subset {\cal A}$  it follows that 
$\partial {\cal A}$ is a non empty union of
unstable leaves and furthermore $S \subset \partial {\cal A}$.
Notice that $S$ is the unique leaf which is either
equal to $G$ or separates $G$ from ${\cal A}$.
In the second case since $g({\cal A}) = {\cal A}$ 
and $g(G) = G$ it follows that $g(S) = S$.
In either case we have that $g(S) = S$.

Then there is an orbit
$\beta$ of $\wwp$ in $S$ with $g(\beta) = \beta$. By theorem
\ref{chain}, $\beta$ and $\gamma$ are connected by a finite 
chain of lozenges $\{ {\cal B}_i \}, 1 \leq i \leq n$. 
Furthermore 

$$E \cap S \not = \emptyset, \ \
E \cap G \not = \emptyset \ \ \Rightarrow E \cap {\cal B}_i
\not = \emptyset, \ \forall i.$$

\noindent
It follows that consecutive lozenges in the chain are adjacent.

\vskip .1in
\noindent
\underline {Claim 2} \ - 
For all $i$, ${\cal B}_i$ is in the front of $\ws(\beta)$.
In particular $\gamma$ is in front of $\ws(\beta)$.

Suppose not. 
Let $r \in \beta$ and $r' \in \gamma$. 
Notice that $p \in \widetilde W^u_+(r')$.
Since $\gamma$ and $\beta$ are
connected by a chain of adjacent lozenges all
intersecting $E$ and $\gamma$ is in the back of
$\ws(\beta)$, it follows that 
$\ooo^u_-(r) = \ooo^u_+(r')$.
For all $i$ big enough $\ws(r)  \cap H_i \not = \emptyset$.
Notice that $g^i(q) \in H_i$.
If $g^i(q)$ is in front of $\ws(r)$ then 
$\ws(g^i(q))$ is in front of $\ws(r)$, contradiction
to $\ws(g^i(q)) = \ws(g^i(p))$ being in the back
of $\ws(r)$. Otherwise $\ws(r) \in \ou(g^i(q))$, implying  
$\ws(r) \in \ou(g^i(p))$ also a contradiction.
This proves claim $2$.


\vskip .1in
Consequently $\gamma$ is in front of $\ws(\beta)$
and $\gamma$, $\beta$ are connected by
and {\underline {even}} number of adjacent lozenges.
Therefore $\ou(r) = \ou(r')$.

Since $R_i$ separates $H_{i}$ from $H_{i-1}$ for all $i$, it follows 
that $\ws(\beta) \cap R_i \not = \emptyset$,
for all $i$ big enough.
Since $g(\ws(\beta)) = \ws(\beta)$  this shows that
$\ws(\beta) \cap H \not = \emptyset$ and similarly
$\ws(\beta) \cap R \not = \emptyset$.
Therefore $\ws(\beta)) \cap {\cal G} \not = \emptyset$
and as a result $\ws(\beta)$ intersects 
$\widetilde W^u_+(\alpha)$.

\vskip .1in
\noindent {\underline {Conclusion:}}
There is an orbit $\beta$ of $\wwp$ with $g(\beta) = \beta$,
$\ws(\beta) \cap \widetilde W^u_+(\alpha)  \not = \emptyset$ and
$\ws(\beta) \cap R \not = \emptyset$.

\blankfig{cros}{2.1}{Impossible intersection of leaves:
a. Case $\delta = \alpha$, 
b. Case $\delta \not = \alpha$.}

Notice that there is $Z \in \fns$ making a perfect fit with
$Y = \wu(\beta)$ so that $Z$ is in the back of $Y$ and
$Z$ and $L$ are not separated, see fig. \ref{cros} a.
Hence $Z, L$ satisfy the hypothesis of the theorem.
As in claim $1$ there is $X \in \fnu$, $X \not = Y$, 
$X$ making a perfect fit 
with $Z$ and intersecting $E$,
see fig. \ref{cros} a.
Therefore the same arguments done before 
work with $G$ replaced by $Y$,
that is the argument works with
$\beta \subset Y$ and $\alpha \subset H$.

Now switch the roles of $Y$ and $H$ and
apply the same argument as above to find
an orbit $\delta$ of $\wwp$ with $h(\delta) = \delta$ and
$\ws(\delta) \cap \widetilde W^u_+(\beta)  \not = \emptyset$,
$\ws(\delta) \cap X \not = \emptyset$.
In addition $\delta$ is connected
to $\alpha$ by an even chain of lozenges
all intersecting a common stable leaf.
Hence if $u \in \delta, u' \in \alpha$, then
$\ou(u) = \ou(u')$.

If $\delta = \alpha$ this produces an
immediate contradiction since  
$\ws(\beta)$ intersects $\widetilde W^s_+(\alpha)$
and $\ws(\alpha)$ intersects $\widetilde W^s_+(\beta)$,
see fig \ref{cros}, a.

Suppose that $\delta \not = \alpha$.
As $\ws(\delta) \cap \widetilde W^u_+(\beta) \not = \emptyset$, 
then $\ws(\beta)$ is 
in the back of $\ws(\delta)$. In particular 
$\ws(\beta) \not \in \ou(u)$.
Hence $\ws(\beta) \not \in \ou(u')$,
a contradiction to the fact that $\ws(\beta)$ 
intersects $\widetilde W^u_+(\alpha)$,
see fig. \ref{cros}, b.

This contradiction implies that $g(H) = H$.
Hence $g(L) = L$ as desired.
\end{proof}

\begin{corollary}{}{}
Let $\Phi$ be an Anosov flow in $M^3$. Suppose
$\Phi$ has branching and $F,L \in \fns$ are not
separated. Then $F$ and $L$ are connected
by an even chain of lozenges, all intersected
by a common stable leaf. In particular
there are only finitely many branching leaves
between $F$ and $L$.
\label{tight}
\end{corollary}

\begin{proof}{}
Up to finite cover we may assume that  $\Phi$ is orientable.
Suppose that $F, L$ are not separated in their negative
sides. 
Let $g \not = id$ be a covering translation
with $g(F) = F$.
By the previous theorem $g(L) = L$.
Let $\gamma$ and $\delta$ be the respective periodic
orbits in $F$ and $L$. 
Furthermore suppose $\wu(\gamma)$ is in the
back of $\wu(\delta)$.

By theorem \ref{chain},
$\gamma$ and $\delta$ are connected by a finite
chain of lozenges. 
Let ${\cal B} = \{ {\cal B}_i \}, 1 \leq i \leq n$,
be the minimal chain from $\delta$ to $\gamma$.
Since $\delta$ is in the back of $\ws(\gamma)$ and
in the front of $\wu(\gamma)$ it follows that
$\gamma$ is the $(+,-)$ corner of ${\cal B}_1$.
Let $\gamma_1$ be the $(-,+)$ corner of ${\cal B}_1$.
Then $\delta$ is in front of $\ws(\gamma_1)$ and in
front of $\wu(\gamma_1)$, hence ${\cal B}_2$ has
$(-,-)$ corner $\gamma_1$ and let $\gamma_2$ be 
the $(+,+)$ corner of ${\cal B}_2$. If $\gamma_2 
= \delta$ we are done. Otherwise  
$\ws(\gamma_2)$ is not separated
from $F$ hence not separated from $L$.
Induction produces $\gamma_4, ... , \gamma_{2k} 
= \delta$ (hence $n = 2k$).
Clearly the $\ws(\gamma_{2i}), 1 \leq i \leq k$
are non separated from each other.

\blankfig{betw}{1.7}{The correct picture of in between branching.}

Conversely suppose that $E \in \fns$ is not separated from
$F, L$ and is between $F$ and $L$.
Let $B_k, k \in {\bf N}$, be a sequence   of stable leaves
so that $B_k \rightarrow F$ as $k \rightarrow \infty$.
As $E$ is not separated from $F$, 
$B_k \rightarrow E$ in $\hs$ when $k \rightarrow \infty$.
But since $F$ and $L$ are connected by a finite chain of lozenges,
then for $k$ big all $B_k$ intersect the interior of 
these lozenges. Therefore
the only possible leaves in the limit of  $B_k$ which are  between $F$ 
and $L$ are those in the stable boundary of the lozenges ${\cal B}_i$.
This completely characterizes such leaves and hence there are
finitely many in between leaves.
\end{proof}

An $\rrrr$-covered Anosov flow can only have
one of two topological types 
(up to isotopy in $\mi$)
for the joint
structure of $\fns, \fnu$ 
\cite{Fe3}. They are characterized by: 

(1) Any leaf of $\fns$ intersects every leaf of
$\fnu$ and vice versa.  This is the called
the {\em product} type.

(2) There is a leaf of $\fns$ which does not intersect every leaf of
$\fnu$. This is the {\em skewed} type, see detailed
definition in \cite{Fe3}.

Suspensions have product type and geodesic flows
have skewed type.

\begin{corollary}{}{}
Let $\Phi$ be an Anosov flow in $M^3$. Then $\Phi$
is topologically conjugate 
to a  suspension of an Anosov diffeomorphism
of the torus if and only if there are no
free homotopies between closed orbits of
$\Phi$ (including non trivial free homotopies
from an orbit to itself).
\label{susp}
\end{corollary}

\begin{proof}{}
If $\Phi$ is not $\rrrr$-covered, theorem \ref{theb} 
shows that there are $F_0 \not = F_1 \in \fns$ and
$g$ a nontrivial covering translation with $g(F_i) = F_i$.
Let $\alpha_i$ be the periodic orbit in $F_i$.
Then $g(\alpha_i) = \alpha_i$. Therefore $\pi(\alpha_0),
\pi(\alpha_1)$ are closed orbits of $\Phi$ (they may be the
same orbit) which are non trivially freely homotopic to
each other.

If $\Phi$ is $\rrrr$-covered and has product
type, then by theorem $2.8$ of \cite{Ba2} (see announcemment in \cite{So})
$\Phi$ is topologically conjugate
to a suspension.
Otherwise $\Phi$ has skewed type and theorem
$3.4$ of \cite{Fe3} produces many non trivial
free homotopies between closed orbits
of $\Phi$.
\end{proof}

Given $2$ adjacent lozenges ${\cal B}_1$ and ${\cal B}_2$
the pivot of their union is the common corner of ${\cal B}_1$
and ${\cal B}_2$.

\begin{corollary}{}{}
Let $\Phi$ be an Anosov flow in $M^3$.
Then up to covering translations there are only
finitely many branching leaves.
\label{fin}
\end{corollary}

\begin{proof}{}
Suppose there are infinitely many inequivalent
stable branching leaves, where the associated
branching is in the positive direction.
Given any two non separated leaves $F, L$
let $\gamma, \alpha$ be the respective 
periodic orbits which are connected by
a chain of lozenges. For any two adjacent lozenges,
the pivot is uniquely determined, furthermore the pivots
are always periodic orbits.

Hence there are infinitely many inequivalent 
periodic pivots $p_i,  i \in {\bf N}$.
Since $\pi(p_i)$ accumlates in $M$, 
assume up to covering translations
that all $p_i$ are in a very small 
product neighborhood of $p \in \mi$, so let 
$i \not = k$ with

$$\wu(p_i) \cap \ws(p_k) \not = \emptyset \ \  {\rm and} \ \
\ws(p_i) \cap \wu(p_k) \not = \emptyset.$$ 

\noindent
An argument exactly like case 1 of theorem \ref{peri}
shows this is impossible.
\end{proof}

We can now completely characterize the structure of the 
set of non separated leaves:

\begin{corollary}{}{}
Let $\Phi$ be an Anosov flow in $M^3$. Let $F$ be a branching 
leaf of $\fns$ and ${\cal E}$ be the set of non separated leaves
from $F$. Given $E, L \in {\cal E}$ we say that
$E < L$ in ${\cal E}$ if there are $G, H \in \fnu$, 
with $G \cap E \not = \emptyset$, $H \cap L \not = \emptyset$ 
and $G$ in the back of $H$. Then either

(1) ${\cal E}$ is finite, hence order isomorphic to $\{1, 2, ... ,n \}$ or,

(2) ${\cal E}$ if infinite and order isomorphic to the set of integers
${\bf Z}$.

\noindent
In particular given any $E, L \in {\cal E}$, there are only finitely
many branching leaves between them.
\label{order}
\end{corollary}

\begin{proof}{}
Up to finite cover if necessary assume that
$\Phi$ is orientable.
Let ${\cal E}$ be the set of non separated leaves from $E \in \fns$.
If ${\cal E}$ is finite, the result is immediate, so assume
it is infinite. Suppose all leaves in ${\cal E}$ are not separated
on their negative sides.


By corollary \ref{fin} there are $E' \not = E^* \in {\cal E}$ 
and $f$ a covering translation with $f(E') = E^*$.
Assume that $E' < E^*$ in the ordering of ${\cal E}$.
Theorem \ref{chain} implies that $E', E^*$ are connected by a
finite chain 
with positive stable boundaries in 
$E_0 = E', E_1, ..., E_n = f(E_0) =  E^*
\in \fns$.
Clearly $E_i < E_j$ if $i < j$.
Since $E_0$ is not
separated from $E_n$, then $f(E_0) = E_n$ is not
separated from $f(E_n)$. This produces 
$E_{n+1},...,E_{2n} = f(E_n)$,
a sequence of non separated leaves. Using $f^i, 
i \in {\bf Z}$, one constructs a sequence $\{ E_i \} _{i \in {\bf Z}}
\subset {\cal E}$
of non separated leaves.

Let now $E \in {\cal E}$. Then $E$ and $E_0$ are not separated, hence
connected by a finite chain of adjacent lozenges all
intersecting a common stable leaf. Notice that the lozenges in the chain
are completely determined by a corner plus a direction.
On the other hand,
starting from $E_0$ and in any direction from $E_0$ (in ${\cal E}$)
there are infinitely 
many adjacent lozenges intersecting a common stable 
leaf.  This implies that $E$ will be 
eventually achieved by lozenges in ${\cal E}$, that is $E = E_i$
for some $i \in {\bf Z}$. 
Hence ${\cal E} = \{ E_i \} _{i \in {\bf Z}}$.
Clearly the order induced above shows that $E_i < E_j$ if $i < j$.
Hence ${\cal E}$ is order isomorphic to ${\bf Z}$ as desired.
\end{proof}

Notice that any covering translation $f$ conjugates
the stabilizers of $F$ and $f(F)$ that is $f \circ ({\cal T}(F)) \circ 
f^{-1} = {\cal T}(f(F))$. Therefore conjugation by $f$ takes a
generator of ${\cal T}(F)$ to a generator of ${\cal T}(f(F))$.

\begin{corollary}{}{}
Let $\Phi$ be an Anosov flow in $M^3$, orientable. 
If ${F_i}, i \in {\bf N}$
$\subset \fns$ is an infinite collection of non
separated leaves of $\fns$, then $M$ has an incompressible
torus.
\label{toro}
\end{corollary}

\begin{proof}{}
As $M$ is orientable, then
if necessary lift to a double cover $M_2$ where
both $\fs$ and $\fu$ are transversely orientable.
The structure  of $\fns, \fnu$ is the same.
By corollary \ref{fin}
there is a covering translation $f$ of $M_2$
with $f(F_i) = F_j$ and $i \not = j$. 

Let $g \not = id$ be 
a generator of the stabilizer of $F_i$
in $\pi_1(M_2)$. 
Then $f g f^{-1}$ is a generator of ${\cal T}(F_j)$.
Theorem \ref{theb}  implies that $g(F_j) = F_j$.
By theorem \ref{indi}, $g$ is indivisible in $\pi_1(M_2)$, 
hence $g$ is also a generator of ${\cal T}(F_j)$.
This implies that either $f g f^{-1} = g$ or $f g f^{-1} = g^{-1}$.

In the first case $f$ and $g$ generate an abelian subgroup
of $\pi_1(M_2)$. If $f^n g^m = 1$, then $f^n g^m (F_i) = F_i$
hence $f^n (F_i) = F_i$. If $n \not = 0$ theorem \ref{indi}
implies that $f (F_i) = F_i$,
contradiction to $F_i \not = F_j$.
Hence $n = 0$. 
Since no multiple of a closed orbit is null homotopic,
them $g^m = id$ implies that $m = 0$ also.
Hence there is a ${\bf Z} \oplus
{\bf Z}$ subgroup of $\pi_1(M_2)$. 

If $f g f^{-1} = g^{-1}$, then
$f^2$ and $g$ generate an abelian subroup of $\pi_1(M_2)$
and the same argument produces ${\bf Z} \oplus {\bf Z} < \pi_1(M_2)$.
Therefore there is a ${\bf Z} \oplus {\bf Z}$ subgroup
of $\pi_1(M)$.
By the torus theorem \cite{Ga} (which uses $M$ being
orientable), either $M$ is a Seifert
fibered space or there is an embedded incompressible
torus. In the first case, Ghys \cite{Gh} proved that
$\Phi$ is up to finite covers, topologically conjugate to
a geodesic flow. But then $\Phi$ would be $\rrrr$-covered, contrary
to hypothesis. Hence $M$ is toroidal as desired.
\end{proof}

\section{Product regions}

\begin{define}{}{}
A positive unstable product region ${\cal P}$
of  $\wwp$ is a region in $\mi$ defined
by an open  strong stable segment 
$\eta \subset F \in \fns$
(or by a flow band $\wwp_{\rrrr} (\eta)$)
so that 

$$ \forall \ p, q \in \eta, \ \ \ooo^u_+(p)  = \ooo^u_+(q).
\ \ \  Then   \ {\cal P} = \bigcup _{p \in \eta} 
\widetilde W^u_+(p).$$

\noindent
The segment \ $\eta$ (which may be infinite)
is called a base segment for the product
region. 
Similarly define negative unstable product regions 
and stable product regions.
\end{define}

The main property of product regions is the following:
for any $F \in \fns$, $G \in \fnu$ so that (i) $F \cap {\cal P}
\not = \emptyset$ and (ii) $G \cap {\cal P} \not = \emptyset$, then
$F \cap G \not = \emptyset$. To see why this is true, notice
first that (ii) implies that $\emptyset
\not = G \cap \eta = p$. By (i) let $q \in \eta$ with 
$F \cap \widetilde W^u_+(q) \not = \emptyset$.
Then $F \in \ou(q)$ hence $F \in \ou(p)$, that is $F \cap 
G \not = \emptyset$. This is the reason for the terminology
product region.


The purpose of this section is to show that the existence
of product regions implies that the flow is $\rrrr$-covered.
The main difficulty is that we will not assume that
$\Phi$ is transitive. With the  additional hypothesis of transitivity
the proof of this fact is simple and was done in \cite{Fe5}.
%
%
%
%
%
%

Given $e > 0$ and $z \in \mi$, let $\tsv(z)$ be the segment
in $\widetilde W^{ss}(z)$ centered at $z$ and with length $e$.

\begin{theorem}{}{}
Let $\Phi$ be an Anosov flow in $M^3$. If there
is a product region in $\mi$
then $\Phi$ is $\rrrr$-covered.
Furthermore any leaf of $\fns$ intersects every leaf of
$\fnu$ and vice versa. As a result $\Phi$ is topologically
conjugate to a suspension Anosov flow.
\label{produ}
\end{theorem}

\begin{proof}{}
By lifting to a finite cover if necessary 
suppose that $\Phi$ is orientable.
Assume that  there is a positive unstable product region 
defined by $\eta \subset \widetilde W^{ss}(y_1)$.
The proof will be achieved by producing bigger
and bigger product regions in $\mi$, which
eventually fill all of $\mi$. This will show
there is a product structure in $\mi$ and hence that 
the flow is $\rrrr$-covered.

If $\Omega$ is the nonwandering set of $\Phi$ then
$W^s(\Omega) = M$ \cite{Pu-Sh}. 
Since the periodic orbits are dense in $\Omega$ \cite{Sm,Pu-Sh}
it follows that the set of annular leaves of $\Phi$ forms a dense
subset of $M$. Therefore there is a periodic orbit
$\gamma$  of $\wwp$ so that if $p \in \gamma$, then
$\widetilde W^u(p) \cap \eta \not = \emptyset$.
If $e > 0$ is small enough then for any $z \in \tsv(p),
\widetilde W^u(z) \cap \eta \not = \emptyset$.
Hence $\tsv(p)$ is the defining segment of a product region.

Let $g$ be  a generator of ${\cal T}(\ws(\gamma))$.
For any $y_2 \in \widetilde W^{ss}(p)$, $y_2$ near
enough $p$, then $y_2 \in \tsv(p)$,
hence $\ooo^u_+(y_2) = \ou(p)$. Since 
$g(\widetilde W^u_+(p)) = \widetilde W^u_+(p)$, then 

$$\ou(g^i(y_2)) = g^i(\ou(y_2)) = g^i(\ou(p)) = \ou(p),
\ \forall i \in {\bf Z},$$

\noindent
Consequently for any $y_3 \in \ws(p)$ it follows that 
$\ou(y_3) = \ou(p)$. 
Let 

$${\cal A} \ \  = \bigcup _{y_2 \in \widetilde W^{ss}(p)} 
\widetilde W^u_+ (y_2).$$

\noindent
Then ${\cal A}$ is a product region with an
infinite basis segment $\widetilde W^{ss}(p)$. 

We now prove that the front of $\ws(p)$ is exactly the set ${\cal A}$.
This shows that there is a product structure 
of $\fns, \fnu$ in the front of $\ws(p)$.

\begin{lemma}{}{}
$\partial {\cal A} = \ws(p)$.
\label{tudo}
\end{lemma}

\begin{proof}{}
Let $a \in \partial {\cal A}$.
Suppose $a \not \in \ws(p)$. There are $a_i \in {\cal A}$
with $a_i \rightarrow a$. Let $b_i = \wu(a_i) \cap 
\widetilde W^{ss}(p)$. Without loss of generality
we may assume $b_i \in \widetilde W^{ss}_+(p)$.
Then $\wu(b_i) \rightarrow \wu(a)$ and maybe
other leaves too.
Notice that $a$ is in front of $\ws(p)$ as all
$a_i$ are.


Suppose that $\wu(a) \cap \ws(p) \not =  \emptyset$.
As $a$ is in front of $\ws(p)$ it would follow 
that $a \in {\cal A}$.
Hence $\wu(a)$ is contained in the front of $\ws(p)$,
in particular $\wu(a) \cap {\cal A} = \emptyset$.

\vskip .1in
\noindent \underline {Claim} - 
$b_i \rightarrow \infty$ in $\widetilde W^{ss}_+(p)$.

Otherwise assume up to subsequence that $b_i 
\rightarrow b_0 \in \widetilde W^{ss}(p)$.
Since $a_i \rightarrow a$ and $\wu(a) 
\cap \widetilde W^{ss}(p) = \emptyset$, 
then $\wu(a), \wu(b_0)$ form a branching pair
of $\fnu$. For $i$ big enough 
$\ws(a) \cap \wu(b_i) \not = \emptyset$.
Hence 

$$\ws(a) \subset \ou(b_i) = \ou(p) = \ou(b_0).$$

\noindent
Hence $\ws(a) \cap \wu(b) \not = \emptyset$, which
is a contradiction to $\wu(a), \wu(b)$ 
being non separated.
The claim follows.

\vskip .1in
Since $\wu(b_i) \rightarrow \wu(a)$;
in fact $\widetilde W^u_+(b_i) \rightarrow \wu(a)$
and also $\wu(a) \cap {\cal A} = \emptyset$ then
$\wu(a) \subset \partial {\cal A}$. Given
$c \in \wu(a)$, choose $c'$ near $c$ with
$c' \in \ws(c)$ and $c' \in {\cal A}$.
Since $\ws(c') \cap {\cal A} \not = \emptyset$ then
$\ws(c') \cap \widetilde W^{uu}_+(p) \not = \emptyset$. As a result
for any  $c \in \partial {\cal A}$ with $c \not \in \ws(p)$, then
$\ws(c) \cap \widetilde W^{uu}_+(p) \not = \emptyset$.
%
%

Let $G = \wu(a)$ and let $F = \ws(p)$.
Notice that $g(G) \not = G$.
Otherwise there is an orbit $\delta$
of $\wwp$ in  $G$ with $g(\delta) = \delta$.
By the above $\ws(\delta) \cap \wu(p) \not = \emptyset$,
a contradiction to both left invariant under $g$.

In fact this shows that $g^n(G) \not = g^m(G)$ for 
any $n \not = m \in {\bf Z}$. Let $G_k = g^k(G)$.
Then $G_k \subset \partial {\cal A}$ so the $G_k$ are 
not separated from each other. 
Therefore by theorem \ref{theb}, $G_k$ contains
a periodic orbit $\delta_k$ and there is an indivisible, 
non trivial covering
translation $f$ with $f(G_k) = G_k$ for all $k \in {\bf Z}$. 
By the above $\emptyset \not = \widetilde W^s_-(\delta_k) 
\cap \widetilde W^{uu}_+(p) = q_k$ for any $k \in {\bf Z}$. Assume
that $q_k = g^k(q_0) \rightarrow p$ as $k \rightarrow +\infty$.

As $\wu(b_i) \rightarrow G_0$, let $S \in \fnu$ with
$S \cap \widetilde W^s_-(\delta_0) \not = \emptyset$
and $S \cap F \not = \emptyset$.
Then $f(S) \cap \widetilde W^s_-(\delta_0) \not = \emptyset$
and we may assume that $f(S)$ is in front of $S$.
As $g$ acts as an expansion in the set of orbits of
$\ws(p)$ then 
$g^j(S) \rightarrow G_0$ as $j \rightarrow +\infty$.
Let $j$ with $g^j(S)$ in front of $f(S)$ and
with $g^j(S) \cap \widetilde W^s_-(\delta_0) \not = \emptyset$,
see fig. \ref{lipro}.
Then $S, g^j(S), \ws(\delta_0)$ and $\ws(p)$ form a rectangle.
As $f(S)$ intersects $\ws(\delta_0)$ between $S \cap \ws(\delta_0)$
and $g^j(S) \cap \ws(\delta_0)$, it follows that  
$f(S) \cap F \not = \emptyset$. In particular
$F$ and $f(F)$ both intersect the unstable leaf $f(S)$.

\blankfig{lipro}{1.8}{Boundaries of product regions.}

If $f(F)$ is in the front of $F$, then as $q_k \rightarrow p$ when
$k \rightarrow +\infty$, it follows that there is some $\ws(q_k)$
which is in the back of $f(F)$, see fig. \ref{lipro}.
This is a contradiction because $f$ leaves
$\ws(q_k)$ invariant. Similarly if $f(F)$ is in the back
of $F$ then $f^{-1}(F)$ intersects $S$ and is in 
front of $F$ producing the same contradiction.

We conclude that $f(F) = F$.
As a result $f = g^n$.
But $f(G_k) = G_k \not = G_{k+n} = g^n(G_k)$,
contradiction.

This shows that the hypothesis $\partial {\cal A}
\not = \ws(p)$ is impossible, hence the lemma follows.
\end{proof}

\noindent
{\underline {Continuation of the proof of theorem
\ref{produ}}} 

%
%

Let $\alpha$ be a periodic orbit, $\alpha \not = \gamma$ with
$\wu(\alpha) \cap \ws(\gamma) \not = \emptyset$.
Let $q \in \alpha$. Assume that $q$ is in front of $\ws(p)$. 
Let ${\cal C} = \bigcup _{z \in \widetilde W^{ss}(q)} 
\widetilde W^u_+(z)$.
Then ${\cal C}$ is a product region and 
as in lemma \ref{tudo}, $\partial {\cal C} = \ws(q)$.
Since $\partial {\cal A} = \ws(p)$,
it follows that ${\cal C} \subset {\cal A}$.

Let $h$ a generator of ${\cal T}(\ws(\alpha))$
so that $h$ acts as an expansion in the set of orbits of $\wu(\alpha)$.
Since $\widetilde W^{uu}_-(q) \cap \partial {\cal A}
\not = \emptyset$,
then for any $i > 0$,
$h^i({\cal A})$ is a product region strictly bigger
than ${\cal A}$ and 
$\partial h^i({\cal A}) = h^i(\ws(p))$.

Therefore for any $z, y \in \widetilde W^{uu}(q)$
there is $i > 0$ so that $z, y \in h^i({\cal A})$.
Let $w = \widetilde W^{uu}(q) \cap h^i(\ws(p))$,
so $z, y \in \widetilde W^{uu}_+(w)$.

If $G \in \fnu$ and $G \in \ooo^s_+(z)$ then
$G$ intersects the front of $\ws(w)$.
By the previous lemma the front of $\ws(w)$ is 
equal to $h^i({\cal A})$. 
As $h^i({\cal A})$ is a product region, then

$$G \cap h^i({\cal A}) \not = \emptyset,
\ \ \ws(y)  \cap h^i({\cal A}) \not = \emptyset \ \ 
\Rightarrow  \ \ 
G \cap \ws(y) \not = \emptyset.$$


\noindent
Since $G$ is in front of $\wu(y)$ then $G \in \ooo^s_+(y)$.
By symmetry  $\ooo^s_+(z) =\ooo^s_+(y)$.

It follows that $\widetilde W^{uu}(q)$ is then
a basis segment of a positive stable product
region ${\cal P}_1$. By lemma \ref{tudo},
$\partial {\cal P}_1 = \wu(q)$.
Similarly $\widetilde W^{uu}(q)$ is also
the basis segment of a negative stable product
region ${\cal P}_2$ and $\partial {\cal P}_2
= \wu(q)$. Hence $\mi = {\cal P}_1 \cup {\cal P}_2$.

It follows from this analysis that for any
$E \in \fns, E \cap \wu(q) \not = \emptyset$.
Therefore $\fs$ is $\rrrr$-covered.
Similarly for any $R \in \fnu$ if it in the 
front of $\wu(q)$, then $R \subset {\cal P}_1$
hence $R \cap \ws(q) \not = \emptyset$,
and similarly for $R$ in the back of $\wu(q)$.
This shows that $\fu$ is also $\rrrr$-covered, hence
that $\Phi$ is $\rrrr$-covered.

Let now $E \in \fns, R \in \fnu$.
Assume that 
$R$ is (say)
in front of $\wu(q)$. 
Then $R \cap {\cal P}_1 \not = \emptyset$ and
$E \cap {\cal P}_1 \not = \emptyset$,
so $E \cap R \not = \emptyset$.
Therefore any leaf of $\fns$ 
intersects every leaf of $\fnu$ and vice versa.
Theorem $2.8$ of \cite{Ba2} implies that $\Phi$ 
is topologically conjugate to a suspension Anosov flow.
\end{proof}

\section{Infinite branching and transverse tori}

In this section we show that, if infinite branching 
occurs, then a particular type of structure, called
a scalloped region, occurs in
$\mi$ (or $\oo$) and there is an embedded torus transverse to
the flow.
We then show that there are many examples
with only finite branching.

\begin{theorem}{}{}
Let $\Phi$ be an Anosov flow in $M^3$. It there is
infinite branching in $\fns$, then there is associated infinite
branching in $\fnu$
\label{infi}
\end{theorem}

\begin{proof}{}
%
Let ${\cal E} = \{ E_i \} _{i \in {\bf Z}}$ be an infinite,
totally ordered collection of non separated leaves.
Assume they are not separated on their negative sides.
Let $\gamma_i$ be the periodic orbit in $E_i$.
Theorem \ref{chain} implies that for any $i$,
$E_i$ forms part of the boundary
of two lozenges:
let ${\cal B}_{2i-1}$ be the lozenge with
$(+,+)$ corner $\gamma_i$ and let ${\cal B}_{2i}$ 
be the lozenge with $(+,-)$ corner $\gamma_i$.
Let $F_i \in \fns$ be the  other leaf in
the boundary of ${\cal B}_{2i}$ and ${\cal B}_{2i+1}$,
where  ${\cal B}_{2i}$ and ${\cal B}_{2i+1}$
are in front of $F_i$.
Let $\zeta_i$ be the periodic orbit in $F_i$, see fig. \ref{forw}.
Then the $\{ F_i \} _{i \in {\bf Z}} \subset \fns$ are all
non separated from each other on their positive sides.
Furthermore all ${\cal B}_i$ intersect 
a common {\underline {stable}} leaf.

\blankfig{forw}{1.7}{Chain of lozenges.}

Notice that the sides of ${\cal B}_{2i}$
are $\widetilde W^s_+(\gamma_i), \widetilde W^u_-(\gamma_i),
\widetilde W^s_-(\zeta_i)$ and $\widetilde W^u_+(\zeta_i)$;
while the sides of ${\cal B}_{2i+1}$
are $\widetilde W^s_-(\gamma_{i+1}), \widetilde W^u_-(\gamma_{i+1}),
\widetilde W^s_+(\zeta_i)$ and $\widetilde W^u_+(\zeta_{i})$,
see fig. \ref{forw}.
Let ${\cal L} = \cup _{i \in {\bf Z}} {\cal B}_i$.
Then all of the following sets are equal:

$$\ooo^u_-(\gamma_i), i \in {\bf Z}, \ \ \ou(\zeta_j), 
j \in {\bf Z}.$$

Let ${\cal C}_i$ be the back of $\wu(\gamma_i)$ and 
let ${\cal C} = \cup _{i \in {\bf N}} {\cal C}_i$.
The set ${\cal L}$ is a union of adjacent lozenges. 
Then for any
$p, q \in \widetilde W^u_-(\gamma_0)$ and any $i > 0$,
$\wu(\gamma_i) \in \ooo^s_+(p) \cap  \ooo^s_+(q)$.
If ${\cal C} = \mi$, then the intersections of $\wu(\gamma_i)$
with $\widetilde W^{ss}_+(p)$
and $\widetilde W^{ss}_+(q)$ are escaping to infinity
in these leaves. This implies that 
$\ooo^s_+(p) = \ooo^s_+(q)$. 
Therefore $\widetilde W^{uu}_-(\gamma_0)$ would be the basis segment
of a positive stable product region in $\mi$. 
By theorem \ref{produ},
$\Phi$ would be $\rrrr$-covered contrary to
hypothesis.
Hence ${\cal C} = \mi$.
This is the key fact which will produce 
a covering translation $f$ commuting with $g$.

Let then $p \in \partial {\cal C}$,
hence $\wu(p) \subset \partial {\cal C}$.
For all $i$ big enough
$\widetilde W^s_-(p) \cap \wu(\gamma_i) 
\not = \emptyset$.
This implies
that $\widetilde W^s_-(p) \cap \widetilde W^u_-(\gamma_i)
\not = \emptyset$ for any $i \in {\bf Z}$.
As a result $\wu(p) \subset \partial {\cal L}$.

Since $g({\cal C}) = {\cal C}$, then
$g^n(\wu(p)) \subset \partial {\cal L}$ for any $n \in {\bf Z}$.
If $g^n(\wu(p)) = \wu(p)$ for some $n \not = 0$,
let $\beta$ be the periodic orbit in $\wu(p)$. Then

$$g^n(\ws(\beta)) = \ws(\beta), \ \  g^n(\wu(\gamma_i)) = \wu(\gamma_i) \ \
{\rm and} \ \ \ws(\beta) \cap \wu(\gamma_i) \not = \emptyset,$$

\noindent 
contradiction.
Hence the leaves
$g^n(\wu(p)), n \in {\bf Z}$ are all distinct and all non separated 
from each other on their negative sides.
By theorem \ref{theb}, $g^n(\wu(p))$ are all periodic and
let $f$  be the indivisible covering translation
leaving all invariant and acting as
an expansion in the set of orbits in $\wu(p)$. 
%

Notice that $g(\ws(p))$ is 
in front of $\ws(p)$. 
Let $H_0 = \wu(p), H_1, ..., H_n = g(\wu(p))$ 
be the chain of non separated leaves from $\wu(p)$
to $g(\wu(p))$.
As in the argument above, 
one constructs $\{  H_k \}_{k \in {\bf Z}}$,
all in $\partial {\cal L}$.
Let $\beta_k$ be the periodic orbits in $H_k$.
Then $\beta_k$ is the corner of two lozenges
${\cal R}_{2k-1}$ and ${\cal R}_{2k}$.
Then all ${\cal R}_k$ intersect a common unstable leaf.

\blankfig{scal}{2.7}{A scalloped region in the universal cover.}

Furthermore if $q \in \partial {\cal C}$, then $\wu(q)$ is not
separated from $H_0$, so $\wu(q)$ is one of $H_k$.
Let $\{ G_k \}_{k \in {\bf Z}}$ be the 
sequence of leaves which form the negative unstable boundary
of the lozenges 
$\{ {\cal R}_k \} _{k \in {\bf Z}}$.
Then $f(G_k) = G_k$ for all $k$.

Given $l \in {\bf Z}$ then for $j > 0$ big enough 
$\widetilde W^u_-(\gamma_j) \cap \widetilde W^s_-(\beta_l)
\not = \emptyset$. Since all $\ooo^s_-(\beta_k), k \in {\bf Z}$
are equal as are all $\ooo^u_-(\gamma_i)$ this implies that for
any $i, k \in {\bf Z}$, ${\cal B}_i \cap {\cal R}_k \not = \emptyset$.
As $g({\cal B}_i) = {\cal B}_i$ for any $i \in {\bf Z}$
and $g({\cal R}_k) = {\cal R}_{k + n}$ for any $k \in {\bf Z}$,
then for any $i \in {\bf Z}$, ${\cal B}_i \subset 
\cup_{k \in {\bf Z}} {\cal R}_k$.

In addition notice that $g^m(\ws(\beta_0)) \rightarrow 
\cup _{i \in {\bf Z}} E_i = {\cal E}$ as
$m \rightarrow +\infty$. As $f(\ws(\beta_k)) = \ws(\beta_k),
\forall k \in {\bf Z}$ then
$f$ leaves invariant the
set ${\cal E}$. Therefore there is $j \in {\bf N}$ so
that $f(E_i) = E_{i + j}$ for all $i \in {\bf Z}$.
Since $f({\cal R}_k) = {\cal R}_k, \forall k \in {\bf Z}$,
then the same argument as above implies that
${\cal R}_k \subset \cup _{i \in {\bf Z}} {\cal B}_i$
for any $k \in {\bf Z}$. We conclude that 

$${\cal L} = \bigcup _{i \in {\bf Z}} {\cal B}_i
= \bigcup _{k \in {\bf Z}} {\cal R}_k.$$

\noindent
The region ${\cal L}$ is called a {\em scalloped} region,
see fig. \ref{scal} and is uniquely associated to the 
infinite branching ${\cal E}$.
Notice that $\fns$ and $\fnu$ restrict to foliations
with $\rrrr$ leaf space in ${\cal L}$.
%
\end{proof}

\begin{theorem}{}{}
Let $\Phi$ be an Anosov flow in $M^3$ orientable. If there is
infinite branching in (say) $\fns$
then there is an embedded torus transverse to $\Phi$.
\label{trans}
\end{theorem}

\begin{proof}{}
Assume first that $\Phi$ is orientable.
We use the notation from the previous theorem.
Let $\nu _{(i,k)} = \widetilde W^u_-(\gamma_i) 
\cap \widetilde W^s_-(\beta_k)$ an orbit of $\wwp$.
Then there are ${\bf Z} \oplus {\bf Z}$ such orbits
in ${\cal L}$.
Recall that 
$g(\gamma_i) = \gamma_i$, $f(\beta_k) = \beta_k$,
and $f$ acts as a contraction in the set of orbits in  
$\ws(\beta_k)$ and likewise for the action $g$  in $\ws(\gamma_i)$.
Then there are $a, b \in {\bf N} - \{ 0 \}$ so that:

$$f(\nu_{(0,0)}) = \nu_{(a,0)}, \ \ {\rm since} \ \
f(\ws(\beta_0)) = \ws(\beta_0),$$
$$g f(\nu_{(0,0)}) = \nu_{(a,b)},  
\ \ {\rm since} \ \ g(\wu(\gamma_a)) = 
\wu(\gamma_a),$$
$$ f^{-1} g f(\nu_{(0,0)}) = \nu_{(0,b)}, \ {\rm as} \ \
f(\ws(\beta_b)) = \ws(\beta_b),
\ f^{-1}(\wu(\gamma_a)) = \wu(\gamma_0),$$
$${\rm and \ \ finally} \ \ 
g^{-1} f^{-1} g f(\nu_{(0,0)})= \nu_{(0,0)}.$$

Since $\nu_{(0,0)} \subset \widetilde W^u_-(\gamma_0)$,
$\nu_{(0,0)}$ is not a periodic orbit of $\wwp$.
Therefore  the last equation above 
implies that $g f = f g$.
Furthermore $f^n g^m = id$, clearly implies that
$n = m = 0$ so $f, g$ generate a ${\bf Z} \oplus {\bf Z}$
subgroup of $\pi_1(M)$. 
Notice that this subgroup preserves ${\cal L}$ and hence
also preserves $\partial {\cal L}$.

Let $p \in \nu _{(0,0)}$. Let $\xi_1$ be an embedded arc 
in $\widetilde W^s_-(\beta_0)$ from
$p \in \nu_{(0,0)}$ to $f(p) \in \nu_{(a,0)}$
transverse to $\wwp$ and so that $\pi(\xi_1)$
is a smooth closed curve in $M$.
Let $\xi _2$ be a similar  arc from
$p \in \nu _{(0,0)}$ to $g(p) \in \nu _{(0,b)}$ contained
is $\widetilde W^u_-(\gamma_0)$.
Since $f g = g f$ then
 $\xi = \xi_1 * f(\xi_2) * 
(g(\xi_1))^{-1} * (\xi_2)^{-1}$
is a closed loop in $\mi$. 
As $\wu(\gamma_0), \wu(\gamma_b), \ws(\beta_0)$
and $\ws(\beta_a)$ form a rectangle it is 
easy to produce an smooth embedded disk $D_1$ in $\mi$,
which is transverse to $\wwp$ and so that 
$\partial D_1 = \xi$.

After a small perturbation of $D_1$ near $\partial D_1$,
we may assume that
$D = \pi(D_1)$ is a smooth closed
surface transverse to $\Phi$.
A priori $D$ is only an immersed surface.
Again after a small perturbation of $D$,
we may assume that $D$ is
transverse to itself. Using cut and paste
techniques
\cite{He,Ja}, as explicit done by Fried \cite{Fr},
one can eliminate all triple points of intersection
and double curves of intersection, transforming
$D$ into a union of embedded surfaces transverse to $\Phi$.

Any such surface has induced stable and unstable foliations
hence it has zero Euler characteristic. It
is transverse to the flow, hence it is two sided in $M$
and as $M$ is orientable, then this
transverse surface has to be a torus.

If $\Phi$ is not orientable, the above proof can
be applied to a double cover $M'$ of
$M$ where the lifted flow is orientable.
The image in $M$ of the transverse torus in $M'$ is an
(immersed) torus in $M$ and again cut and
paste techniques  yield the result.
\end{proof}

\vskip .1in
\noindent
{\bf Remarks:} 
(1) With much more work, one can in fact show there is a transverse 
torus $T$ intersecting exactly those orbits in $\pi({\cal L})$.
This is done in detail by Barbot in
the proof of theorem $B$ in \cite{Ba3}
where the hypothesis are the existence of 
two commuting covering translations $f, g$,
so that both $f$ and $g$ are associated to (different
free homotopy classes of) closed
orbits of $\Phi$ in $M$.



(2) As mentioned earlier the Anosov flow constructed by
Bonatti and Langevin has infinite branching.
The scalloped  region of this flow was explained  
in  detail in \cite{Fe5}.
The Bonatti-Langevin flow is the 
simplest Anosov flow with infinite branching
in the sense that there is only one orbit $\gamma$
of $\Phi$ which does not intersect the transverse torus $T$. 
In  this case all periodic orbits in $\partial {\cal L}$
are  lifts of $\gamma$.
The picture in $\mi$ is very symmetric.

\vskip .1in
We now prove that a large class of Anosov flows
in dimension $3$ have only finite branching.
If the branching is finite we define its length
to be the number of non  separated leaves.

\begin{theorem}{}{}
Any Anosov flow obtained by the Franks-Williams construction
\cite{Fr-Wi} is not $\rrrr$-covered and has only
length $2$ branching.
\end{theorem}

\begin{proof}{}
First we recall the Franks-Williams construction \cite{Fr-Wi}.
Start with a suspension Anosov flow $\Phi_0$ in $N$ and a 
closed orbit $\overline \gamma$.
Modify the flow in a neighborhood of $\overline
\gamma$ using Smale's DA (derived from Anosov) construction
\cite{Sm,Wi},
so that $\overline \gamma$ becomes an expanding
orbit and $2$ new hyperbolic orbits $\gamma_1$ and
$\gamma_2$ parallel to $\overline \gamma$ are created, see
fig. \ref{der}, a.
This produces the new flow $\Phi^*$ in $N$.

One can do this in a way that the stable foliation
is still preserved by the new flow. 
Now remove a solid torus neighborhood $V$ of $\overline \gamma$
with a boundary torus $T_1$, transverse to the flow.
This creates a manifold $M_1 = N - V$ where the flow
is incoming in the boundary. There is an induced
stable foliation in $\partial M_1 = T_1$,
which has two closed leaves and two Reeb components
in between them.
Using a time reversal of this flow
construct $M_2$ with a boundary torus
$T_2$ where the flow is outgoing and there
is an induced unstable foliation in $T_2$.
Finally glue $T_1$ to $T_2$ so that after glueing
the stable and unstable foliations are
transverse.
Let $T$ be the torus obtained by glueing $T_1$ to $T_2$.
Franks and Williams show that such flows
are Anosov and clearly intransitive since
$T$ is a separating torus.
Hence the flows are not $\rrrr$-covered.

By theorem \ref{theb}, any branching of $\fns$ and
$\fnu$ produces freely homotopic closed orbits,
so we first understand free homotopies.
Let $\alpha$ and $\beta$ be freely homotopic closed orbits
of $\Phi$. 
and let $\tau: A \rightarrow M$ be an annulus realizing
the free homotopy.
Assume that $A$ is in general position 
and is transverse to $T$. Notice that
$\partial A = \alpha \cup \beta$ is disjoint from $T$.
Then $\tau ^{-1}(T)$ is a union of closed curves in $A$.
We can eliminate all null homotopic components
as follows: since $T$ is transverse to $\Phi$, it
is incompressible \cite{Fe4}. Then any null homotopic component
of $\tau^{-1}(T)$ also produces a null homotopic curve in $T$.
Using cut and paste arguments \cite{He,Ja} and the fact that 
$M$ is irreducible we can eliminate this component by a homotopy
of the annulus.
We may then assume that $\tau^{-1}(T)$ is a union of 
finitely many curves parallel to $\partial A$.

\blankfig{der}{2.2}{a. DA construction, b. Induced foliations
in a lift of the torus.}

Let now $B_1$ be the closure of a component 
of $A - (\tau^{-1}(T))$  containing
a boundary component of $A$.
Let $\partial _1 B$ be this boundary component 
(suppose that $\tau(\partial _1 B) = \alpha$)
and let $\partial _2 B = \partial B - \partial _1 B$.
Assume that $\tau(B_1) \subset M_1$.
Notice that $M_1$ fibers over the circle
with fiber $F$ a torus minus a disk.
Any closed orbit of $\Phi$ in $M_1$ has non zero algebraic
intersection with $F$, hence the same is true
for the other boundary of $B$, that is $\partial _2 B$ is
not a multiple of the meridian.
Reglue the solid torus $V$ as originally
to recover $N$ and the DA flow $\Phi^*$ in $N$.
Since $\partial _2 B$ is not a meridian, 
then $\partial _2 B$ is freely
homotopic (in $V$) to $\overline \gamma^n, n \not = 0$, hence
freely homotopic to $\gamma^n_1$.

The DA construction is equivalent to
spliting $\wu(\overline \gamma)$ into 
two and blowing air in between the 2 sides \cite{Wi}, much 
in the sense of essential laminations \cite{Ga-Oe}.
In particular there is a topological semiconjugacy 
between $\Phi^*$ and $\Phi_0$.
Hence free homotopies between closed orbits of
$\Phi^*$ produce a free homotopy between  two closed orbits
of $\Phi_0$ in $N$.
But any free homotopy in a suspension is trivial \cite{Fe3}.
Therefore $\alpha$ is either $\gamma_1$ or $\gamma_2$
and the free homotopy can be homotoped into $W^s(\gamma_1)$
(or into $W^s(\gamma_2)$).

Furthermore $M_1$ is acylindrical,
that is, any properly immersed annulus can
be homotoped into the boundary.
This is due to Waldhausen (for a proof see \cite{Jo})
and follows from the fact that $M_1$ is atoroidal
(in fact it is hyperbolic \cite{Th2}), $M_1$ not a Seifert fibered
space and $\partial M_1$
is a single torus.

These two facts 
imply that the only non trivial
free homotopies between closed orbits of
$\Phi$ can always be homotoped into $T$.
Notice that there are such free homotopies,
since $\gamma_1$ is freely homotopic to $\gamma_2$
in $M_1$ and also there are two closed 
orbits  of $\Phi$ in $M_2$ which are freely homotopic
to each other and freely homotopic to $\gamma_1$. 
These orbits are associated to the 4 closed leaves of 
the induced stable and unstable foliations in $T$.

As a result of this, in order to  understand branching in the 
universal cover all we need to do is understand
the structure of $\fns,\fnu$ induced in 
lifts of $T$.
Since there are two closed leaves in $\fs \cap T$
and Reeb components in between them 
and similarly for $\fu \cap T$,
the picture in the universal cover of
$T$ is as in fig. \ref{der} b.
This shows that $F_1 \in \fns$ is
not separated from $F_2$ on their negative sides,
$F_2$ not separated from $F_3$ in their positive
sides and so on. 
This implies that  $F_1$ is separated from
$F_3$.
Therefore $F_2$ is the only leaf non separated from
$F_1$ which is in the negative side of $F_1$.
We conclude that in such flows any branching has 
length two.
\end{proof}

\section{Continuous extension of Anosov foliations}

If $\Phi$ is an Anosov flow in $M^3$, 
Sullivan \cite{Su} showed  that the intrinsic geometry
of leaves of $\fns$ and $\fnu$ is 
{\em negatively curved in the large}
as defined by Gromov \cite{Gr}.
This holds without any assumption on $M$.
Then any leaf $F \in \fns \cup \fnu$ has a canonical
compactification with an intrinsic ideal boundary 
$\pin F$ \cite{Gr}.
We proved in \cite{Fe2} that $\pin F$ is always
homeomorphic to a circle. 

If $F \in \fns$ then the intrinsic ideal points
correspond to the (distinct) negative limit points
of flow lines in $F$ and to the common
positive limit point of all flow lines \cite{Fe3}.
The intrinsic geometry of $F \in \fns$ resembles
that of the hyperbolic plane ${\bf H}^2$ where the flow
lines correspond to the geodesics in ${\bf H}^2$
which have a common limit point in the ideal
boundary of ${\bf H}^2$,
see fig. \ref{intr}.
Analogous results hold for $\fnu$.

\blankfig{intr}{1.8}{Intrinsic ideal points.}

If $p \in F  \in \fns$, we define 
$p_- \in \pin F$ to be the intrinsic negative limit point
of the flow line through $p$, that is
$p_- = \lim_{t \rightarrow -\infty} \wwp_t(p)$,
where the limit is taken in $F \cup \pin F$,
see fig. \ref{intr}.
Similarly define $p_+$.
For any $p, q \in F \in \fns$, $p_+ = q_+ \in \pin F$ 
and this is also denoted by $F_+$.
Furthermore if $p_i \in \widetilde W^{ss}(p)$
and $p_i \rightarrow \infty$ in $\widetilde W^{ss}(p)$,
then $(p_i)_- \rightarrow p_+$ as points in $\pin F$
\cite{Fe3}. This can be clearly seen
in the model of the hyperbolic plane.

Notice that $p_- = q_-$ for any $p, q$ in the same flow
line $\alpha$ of $\wwp$, so this is also denoted by
$(\alpha)_- \in \pin F$ and similarly $(\alpha)_+ = F_+$.

From now on we assume that 
$\pi_1(M^3)$ is
{\em negatively curved} as defined by Gromov \cite{Gr}.
Gromov constructed a canonical compactification
of $\mi$ with an ideal boundary $\partial \mi$.
When $M$ is irreducible (always the case for us),
Bestvina and Mess \cite{Be-Me} showed that $\partial
\mi$ is homeomorphic to a sphere, denoted by $\si$.
Furthermore $\mi \cup \si$ is homeomorphic to
a closed $3$-ball.

Recall that the foliations $\fns, \fnu$ 
are transversely oriented.

\begin{define}{}{}
The limit set of a subset $B$ of
$\mi$ is $\Lambda _B = \overline B \cap \si$,
where the closure is taken in $\mi \cup \si$.
Given $F \in \fns$ or $\fnu$ 
and $p \in \si - \Lambda_F$,
we say that $p$ is above $F$ if there
is a neighborhood $U$ of $p$ in $\mi \cup \si$
so that $U \cap \mi$ is in front of $F$.
Otherwise we say that $p$ is below $F$.
Given a connected component of $\si - \Lambda_F$ either
all of its points are above $F$ and we say
this component is above $F$ or all points
are below $F$ and we say the component is below $F$.
Similarly for $G \in \fnu$.
\end{define}

%
%

\begin{proposition}{}{}
Let $\Phi$ be an Anosov flow in $M^3$ with negatively
curved $\pi_1(M)$. Either $\Lambda_F = \si$
for every $F \in \fns$;
or for every $F \in \fns$,
$\si - \Lambda_F$ has at least one connected component above $F$
and one component below $F$.
\label{compo}
\end{proposition}

\begin{proof}{}
Classical $3$-dimensional topology \cite{He,Ja} and Smale's
spectral decomposition theorem \cite{Sm} imply that $\Phi$
is transitive \cite{Fe4}.

We may assume that $\fs, \fu$ are transversely orientable.
Suppose there is $F \in \fns$ with $\Lambda_F \not
= \emptyset$.
Assume that there is a component $Z$ of $\si - \Lambda_F$
which is above $F$.
Since stable leaves are dense in $M$, 
then for every $L \in \fns$ there is a covering
translation $g$ with $g(L)$ in the back
of $F$ and so that $F, g(L)$ intersect a common unstable
leaf. Then $Z \cap \Lambda_{g(L)} = \emptyset$
and since $Z$ is in front of $F$, it is also
in front of $g(L)$. Therefore there is a component
of $\si - \Lambda_{g(L)}$ above $g(L)$. Translating
by $g^{-1}$ we conclude that there is a component
of $\si - \Lambda_L$ above $L$.

If $\Phi$ were $\rrrr$-covered, 
then $\Lambda_F = \si$ \cite{Fe2}, contrary to assumption.
Hence $\Phi$ is not $\rrrr$-covered  and 
by transitivity, it follows that
$\Phi$ has branching in the positive and negative
directions \cite{Fe5}. Let then $E, E' \in \fns$ so that
they are not separated on their negative sides.
By the above argument $\si - \Lambda_E$ has a component
$Z_0$ above $E$. Since $E$ is in the back of $E'$ and
$E'$ is in the back of $E$ it follows that
$Z_0 \cap \Lambda_{E'} = \emptyset$ and all points
in $Z_0$ are below $E'$. Hence $\si - \Lambda_{E'}$
has a component below $E'$. Using the same argument
as above we conclude that for every $L \in \fns$,
there is a component of $\si - \Lambda_L$ below $L$.
This completes the proof.
\end{proof}

We say that $\wwp$ has the
{\em continuous extension property}
if for any leaf $F \in \fns \cup \fnu$,
the embedding $\varphi _F: F \rightarrow \mi$,
extends continuously to $\varphi _F: F \cup \pin F
\rightarrow \mi \cup \si$.
This gives a continuous parametrization of
the limit sets $\Lambda_F = \varphi _F (\pin F)$.
This also implies that there is a continuous function

$$\eta_-: \mi \rightarrow \si, \ \ \ 
\eta_-(x) = \lim_{t \rightarrow -\infty} \wwp_t(x),$$

\noindent
where the limit is computed in $\mi \cup \si$.
Since the function is constant along an orbit
$\alpha$ of $\wwp$, this will also denote
$\et(\alpha)$.
Furthermore for any $G \in \fnu$, $\eta_-$ is
a constant function in $G$ with value $\varphi _G (G_-)$.
Similarly define $\eta_+: \mi \rightarrow \si$.
The continuous extension property implies that for any
$p \in F \in \fns$,
$\Lambda _F =  \varphi _F (\pin F) =
\et(\widetilde W^{ss}(p)) \cup \eta_+(p)$.

In \cite{Fe6} we study 
the continuous extension property for $\rrrr$-covered flows.

\begin{theorem}{}{}
Let $\Phi$ be an Anosov flow in $M^3$ with
negatively curved $\pi_1(M)$.
Suppose that $\Phi$ is not $\rrrr$-covered
and in addition that $\wwp$ has the continuous
extension property. Then for any leaf 
$C \in \fns \cup \fnu$, the limit set $\Lambda _C$
is a Sierpinski curve, that is the complement of a
countable, dense union of open disks in the sphere $\si$.
\label{sier}
\end{theorem}

\begin{proof}{}
We may assume that $\fs, \fu$ are transversely orientable.
We first prove that $\Lambda_C \not = \si$ and then use part of
the proof of this fact to show that limit sets are Sierpinski curves.
The first part is the same as the proof of theorem $5.5$ of \cite{Fe4}.
In \cite{Fe4} we used the hypothesis of quasigeodesic behavior 
of flow lines of $\wwp$ in order to describe the structure
of branching of $\fns$ and $\fnu$. In this article we
obtained a description of branching without any hypothesis
and this is what is needed for the proof of theorem \ref{sier}.

Since $\Phi$ is transitive, $\fns$ has branching in
the positive and negative directions.
Using theorem \ref{order} we produce
$\Theta$, a union of two adjacent lozenges
in $\mi$ (or ${\cal O}$) intersecting a common stable
leaf  so that: (1)  the boundary of $\Theta$ has
unstable sides in $G, S \in \fnu$, 
and stable sides in $E, F, L \in \fns$
(2) $E, L$ are not separated on their negative sides,
(3) $G$ is in the back of $S$ and (4) $E \cap G \not = \emptyset$,
$L \cap S \not = \emptyset$,
see fig. \ref{limits}.
By $G$ we mean the half leaf in the boundary
of $\Theta$. Then $\pi(G)$ is dense in $M$ \cite{Fe2}.

\blankfig{limits}{2.5}{Sequence of lozenges.}

Let $C \in \fns$ be a leaf intersecting both $G$ and
$S$, hence $C$ intersects $\Theta$.
Choose a  covering translation $g_1$
so that 

$$g_1(G) \cap F \not = \emptyset, \ \ g_1(G) \cap L \not
= \emptyset.$$

\noindent
Since $g_1(F)$ makes a perfect fit with $g_1(G)$,
then $g_1(F)$ is in the back of $F$.
Since $g_1(L)$ makes a perfect fit with $g_1(E)$ then both
are in the front of $L$. 
Finally
$g_1(S)$ is in the front of $g_1(G)$, in the back
of $S$ and intersects both $L$ and $F$.
Inductively choose covering
translations $g_i$ so that 
$g_i(G)$ is in the back of $S$,

$$g_i(G) \cap F \not = \emptyset, \ \ g_i(G) \cap L \not = 
\emptyset, \ \ g_i(G) \rightarrow S \ 
{\rm as} \ i \rightarrow \infty,$$

\noindent
and $g_i(G)$ is in the front of $g_{i-1}(S)$, see fig.
\ref{limits}.
Let $G_i = g_i(G)$ and similarly define $F_i, L_i, S_i$ and $E_i$.

Let $C_i = C  \cap  g_i(\Theta)$.
\ \ For any flow line $\gamma \in F_i$,
$\wu (\gamma)$ intersects $C_i$ and vice versa.
Hence $\eta_-(C_i) = \eta_-(F_i)$.
Let $q \in C \cap S$.
By continuity of $\eta _-$, there is a
neighborhood $Y$ of $q$ in $\mi$
so that $\eta _-(Y)$ is contained in
a small neighborhood $Y'$ of $\eta _-(q)$ in $\si$.
As $C_i \cap \widetilde W^{ss}(q) 
\rightarrow q$, then
$\eta_-(C_i) \subset Y'$ for $i$ big enough.
Therefore $\eta_-(F_i) \subset Y'$
and as a result $\Lambda_{F_i}$ is 
contained in the closure of $Y'$ and
is not $\si$.

We can now apply the previous proposition to deduce that
for any $L' \in \fns$,
there are components of $\si - \Lambda_{L'}$
above $L'$ and components below $L'$.

For each  $i$ let
$Z_i$ be a component of $\si - \Lambda_{F_i}$
below $F_i$. Since $C$ is in front of $F_i$,
$Z_i \cap \Lambda_C = \emptyset$.
Hence $Z_i$ is contained in a component $Z^*_i$ of
$\si - \Lambda_C$ which is below $C$.
The argument above used to prove that $\Lambda_C \not = \si$
shows that $\Lambda_{F_i}
\subset \Lambda_C$, hence the component $Z^*_i$ of
$\si - \Lambda_C$ is equal to $Z_i$.


For each $i$, $Z_i$ is below $F_i$.
In addition for each  $i \not = j$, 
$F_i$ is in the front of $F_j$ and $F_j$
is in the front $F_i$.
This implies that $Z_i \cap Z_j = \emptyset$.
Hence $\{ Z_i \}, i \in {\bf N}$ is an
infinite family of distinct components of $\si - \Lambda_C$
below $C$.
Using branching of $\fns$ in the negative direction,
one constructs countably 
many components of $\si - \Lambda_C$ above $C$.

Since $\Phi$ is transitive, then
for any $C' \in \fns$ there is a covering translation
$f$ so that $f(C') \cap \Theta \not = \emptyset$.
Since $\si - \Lambda_{f(C')}$ has infinitely many components
above and below $f(C')$, translation by $f^{-1}$ yields
the same result for $C'$.

Notice that these arguments also imply that for every leaf
$F' \in \fns$ either $F'$ is in the back of
$C$, hence $\Lambda_{F'}$ misses at least all
components of $\si - \Lambda_C$ above $C$ for fixed $C$;
or $F'$ is in front of $C$ and $\Lambda_{F'}$ misses
all components of $\si - \Lambda_C$ below $C$ for a fixed $C$.
In particular this implies that there is $\epsilon > 0$ so
that every $\Lambda_{F'}$ misses at least {\underline {some}} disk
of radius $\epsilon$ in $\si$.

Suppose now that for some $R$ in $\fns$,
$\Lambda_R$ has no empty interior.
Let $h$ be a covering translation with both fixed
point in the interior of $\Lambda_R$. By
applying $h^n$ for $n$ big we get $\Lambda_{g^n(R)}$
is almost all of $\si$ except for an 
arbitrarily small neighborhood of the attracting
fixed point of $h$. This contradicts the previous
paragraph. This finishes the proof of the theorem.
\end{proof}



\begin{lemma}{}{}
Let $\Phi$ be an Anosov flow in $M^3$, with
$\pi_1(M)$ negatively curved.
Suppose that $\wwp$ has the continuous extension 
property.
If $F \in \fns$ is periodic and $x \in F$ is
in the periodic orbit of $F$ then for any
$x_1 \in \widetilde W^{ss}(x)$,
$\eta_-(x_1)  = \eta_-(x)$ implies that
$x_1 = x$.
Furthermore $\eta_-(x_1) \not = \eta_+(x)$.
\label{self}
\end{lemma}

\begin{proof}{}
Let $h$ be the generator of ${\cal T}(F)$
associated to  the closed orbit
$\pi(\wwp_{\rrrr}(x))$ traversed 
in the positive flow direction.
Hence $h$ acts as an expansion in the set of
orbits of $\wwp$ in $F$.
Suppose that $\eta_-(x_1) = \eta_-(x)$, but $x_1 \not \in
\gamma = \wwp_{\rrrr}(x)$.
Let $\alpha = \wwp_{\rrrr}(x_1)$.
Then

$$\eta_-(\gamma) = \eta_-(\alpha) \ \
\Rightarrow \ \ 
\eta_-(h^n(\alpha)) = h^n(\eta_-(\alpha)) = 
h^n(\eta_-(\gamma)) = \eta_-(\gamma).$$

\noindent
But 

$${\rm But} \ \ \lim_{n \rightarrow +\infty} 
\eta_-(h^n(\alpha)) \ = \
\lim _{n \rightarrow +\infty} \varphi ( (h^n(\alpha))_- ) \ = \
\varphi (F_+) \ = \
\eta_+(\gamma),$$

\noindent
because the intrinsic negative
limit points of $h^n(\alpha)$ converge in $\pin  F$
to the positive limit point associated to $F$
and in addition $F$ extends continuously to $\si$.
This would imply  $\eta_-(x) = \eta_+(x)$, a contradiction
to $x$ being in a periodic orbit.
Similarly $\eta_-(x_1) \not = \eta_+(x)$.
This proves the lemma.
\end{proof}

We now prove a local property of the limit sets.
Given $L \in \fns$, $\Lambda_L$ is the image
of $\partial _{\infty} L \simeq {\bf S}^1$
under a continuous map, hence $\Lambda_L$ is
locally connected.

\begin{theorem}{}{}
Let $\Phi$ be an Anosov flow in $M^3$ with
negatively curved fundamental group.
Assume that $\Phi$ has the continuous extension
property.
Given any $L \in \fns$ or $\fnu$ and 
any $p \in \Lambda_L$, then for each neighborhood
$U$ of $p \in \si$, $\Lambda_L \cap U$ is neither 
a Jordan arc nor a Jordan curve.
\end{theorem}

\begin{proof}{}
Since $\Lambda_L$
is locally connected, we can choose $U$ so
that $U \cap \Lambda_L$ is connected.

Assume that $U \cap \Lambda_L$ is a Jordan arc or Jordan curve
and let $z$ be a relative  interior point. Suppose that 
$z \not = \eta_+(L)$.
Then $z = \eta_-(c')$ for some $c' \in L$
and there is $e' > 0$ small enough so that
the segment $\sigma ^s_{e'} (c') \subset
\widetilde W^{ss}(c')$ 
(as defined in section $5$), satisfies
$\eta_-(\sigma ^s_{e'}(c')) \subset U \cap \Lambda_L$.
As this is connected we can assume this is a Jordan arc.
Since periodic orbits are dense in $M$,
choose $c$ in a periodic orbit of $\wwp$ near
$c'$ and let $e > 0$ so that the segment
$\sigma ^s_e(c) \subset \widetilde W^{ss}(c)$ satisfies
the following property:
$G \in \fnu$ intersects $\sigma^s_{e'}(c')$ if and only
if it intersects $\sigma^s_e(c)$. Then $K = \eta_-(\sigma^s_{e'}(c'))
= \eta_-(\sigma^s_e(c))$ is a Jordan arc.
We may also assume that $\eta_+(c) \not \in K$.
Let $C = \ws(c)$.

%
%
%
%

Let $g$ be the covering translation associated to the closed 
orbit $\pi(\wwp_{\rrrr}(c))$ of $\Phi$ and assume that
$\eta_-(c)$ is the repelling fixed point of $g$.
Therefore $g$ does not fix the endpoints of
$K$ and $g(K)$ is a Jordan arc strictly bigger
than $K$ in both directions.
Notice that $g(K) \subset \Lambda_C$.

Then 
$g^i(K) =
\et(g^i(\tsv(c)))$ 
is also a Jordan arc $\forall i \in {\bf N}$.
Express $g^i(K)$ as the image of an embedding
$\tau_i : [-i,i] \rightarrow \si$, so
that if $x \in  [-i,i]$ and $j > i$, then
$\tau_i(x) = \tau_j(x)$.
Since $g$ acts as an expansion in the set
of orbits of $\wwp$ in $C$ then

$$\forall \ e^* > 0,  \ \  \exists \
i > 0 \ \ | \ \   \et(\sigma^s_{e^*}(c))
\subset g^i(\et(\tsv(c))) = g^i(K) \ \ \ (*).$$

\noindent
Therefore we can express 
$\et(\widetilde W^{ss}(c))$ as the 
image of an embedding $\tau: {\bf R} = 
{\bf S}^1 - \{ \infty \} \rightarrow \si$,
so that if $x \in [-i,i]$ then $\tau(x) = \tau_i(x)$.

Define $\tau(\infty) = \eta_+(c)$.
By lemma \ref{self}, $\tau:{\bf S}^1 \rightarrow \si$
is injective.
Clearly $\tau$ is continuous in ${\bf S}^1
- \{ \infty \}$.
Suppose that $x_i \in {\bf S}^1$ and $x_i
\rightarrow \infty$, $x_i \not = \infty$.
Choose $c_i \in \widetilde W^{ss}(c)$
with $\tau(x_i) = \et(c_i)$.
If there is a subsequence $i_k$ where $c_{i_k}$ is
bounded in $\widetilde W^{ss}(c)$
then assume this is the  original sequence.
Then by $(*)$ there is $j_0  \in {\bf N}$
so that $\et(c_i) \subset g^{j_0}(K)$
for all $i$ in ${\bf N}$.
Since $\tau$ is injective this implies
that $x_i \in [-j_0,j_0]$ contradiction
to $x_i \rightarrow \infty$.

Hence $c_i \rightarrow \infty$ in 
$\widetilde W^{ss}(c)$ and as $C$
extends continuously to $\si$

$$\tau(x_i) = \et(c_i) \ \rightarrow 
\ \eta_+(c) = \tau(\infty) \ \ {\rm as} \ \ 
i \rightarrow \infty.$$

\noindent
Therefore $\tau$ is continuous, hence
a homeomorphism.
This implies that $\tau({\bf S}^1)
= \Lambda _C$ is a Jordan curve.
If $\Phi$ is $\rrrr$-covered this contradicts
the fact that $\Lambda_C = \si$. If
$\Phi$ is not $\rrrr$-covered this contradicts
theorem \ref{sier}.
The result follows.
\end{proof}

\section{Non $\rrrr$-covered Anosov flows
in hyperbolic $3$-manifolds}

\begin{theorem}{}{}
There is a large class of non $\rrrr$-covered
Anosov flows in hyperbolic $3$-manifolds,
including all Anosov flows in non orientable
hyperbolic $3$-manifolds.
\label{hyper}
\end{theorem}

\begin{proof}{}
Theorem $C$ of \cite{Ba2} states that if 
$\Phi$ is an $\rrrr$-covered Anosov flow
in $M^3$, then either $\Phi$ is topologically
conjugate to a suspension Anosov flow or
the underlying manifold is orientable
(notice that Barbot uses the term ``product"
instead of $\rrrr$-covered).
Since hyperbolic manifolds can never be
the underlying manifolds of suspension
Anosov flows, it suffices to
produce Anosov flows in non orientable
hyperbolic $3$-manifolds.

Consider therefore the suspension of an
orientation  reversing Anosov diffeomorphism
of the torus $T^2$.
Let $M$ be the underlying manifold of the suspension
and let $\alpha$ be an orientation preserving
closed orbit of the flow.
As described by Goodman \cite{Go} and Fried \cite{Fr},
one can do Dehn surgery along this orbit.
Then $(n, 1)$ Dehn surgery on $\alpha$ yields
an Anosov flow in the surgered manifold
$M_{(n, 1)}$. 

Notice now that $(M - \alpha)$ is irreducible, 
atoroidal and homeomorphic to the interior
of a compact $3$-manifold with boundary.
By Thurston's hyperbolization theorem \cite{Th2,Mor}
it follows that $(M - \alpha)$ admits a 
complete hyperbolic structure of finite 
volume. By the hyperbolic Dehn surgery
theorem \cite{Th1}, most Dehn fillings on $(M - \alpha)$
yield closed, hyperbolic manifolds. Since $M$
was non orientable, all of these manifolds
are non orientable.
Whenever the Dehn surgery coefficient is of the form
$(n,1)$, the surgered manifold admits an Anosov flow.
This produces infinitely many Anosov flows in
non orientable hyperbolic $3$-manifolds
and finishes the proof.
\end{proof}

{\small
{

\setlength{\baselineskip}{0.4cm}

}
}

\end{document}